\pdfoutput=1
\documentclass[12pt]{article}
\usepackage[margin=28mm]{geometry}
\usepackage{indentfirst}
\usepackage[usenames,dvipsnames,svgnames,table]{xcolor}
\usepackage[T1]{fontenc}
\usepackage{amsfonts,amsmath,amssymb,amsbsy,amsthm}
\usepackage{hyperref}
\hypersetup{ 
unicode=true,
colorlinks=true,
linkcolor={black},
citecolor={black},
urlcolor={blue!60!black},
pdftitle={Clustered Variants of Haj\'os' Conjecture},
pdfauthor={Chun-Hung Liu and David~R.~Wood}}
\usepackage{enumitem}
\setlist[itemize]{itemsep=0ex,parsep=0.5ex}
\setlist[enumerate]{itemsep=0ex,parsep=0.5ex}
\setlist[description]{itemsep=0ex,parsep=0.5ex}
\usepackage[longnamesfirst,numbers,sort&compress]{natbib}
\makeatletter
\def\NAT@spacechar{~}
\makeatother
\usepackage[noabbrev,capitalise]{cleveref}
\crefname{lem}{Lemma}{Lemmas}
\crefname{thm}{Theorem}{Theorems}
\crefname{prop}{Proposition}{Propositions}
\crefname{conj}{Conjecture}{Conjectures}
\crefformat{equation}{(#2#1#3)}
\Crefformat{equation}{Equation #2(#1)#3}
\newtheorem{theorem}{Theorem}
\newtheorem{lemma}[theorem]{Lemma}

\newtheorem{corollary}[theorem]{Corollary}
\newtheorem{claim}{Claim}[theorem]
\newcommand{\ngs}[2][]{N_{#1}^{\geq s}(#2)}
\newcommand{\nls}[2][]{N_{#1}^{< s}(#2)}
\newcommand{\nlss}[2][]{N_{#1}^{< \beta'}(#2)}
\makeatletter
\renewcommand\section{\@startsection {section}{1}{\z@}{-3ex \@plus -1ex \@minus -.2ex}{2ex \@plus.2ex}{\normalfont\large\bfseries}}
\renewcommand\subsection{\@startsection{subsection}{2}{\z@}{-2.5ex\@plus -1ex \@minus -.2ex}{1.5ex \@plus .2ex}{\normalfont\normalsize\bfseries}}
\renewcommand\subsubsection{\@startsection{subsubsection}{3}{\z@}{-2ex\@plus -1ex \@minus -.2ex}{1ex \@plus .2ex}{\normalfont\normalsize\bfseries}}
 \renewcommand\paragraph{\@startsection{paragraph}{4}{\z@}{1.5ex \@plus.5ex \@minus.2ex}{-1em}{\normalfont\normalsize\bfseries}}
\renewcommand\subparagraph{\@startsection{subparagraph}{5}{\parindent}  {1.5ex \@plus.5ex \@minus .2ex}  {-1em} {\normalfont\normalsize\bfseries}}
\makeatother
\def\C {{\mathcal C}}
\def\T {{\mathcal T}}
\def\G {{\mathcal{G}}}
\renewcommand{\geq}{\geqslant}
\renewcommand{\leq}{\leqslant}
\newcommand\abs[1]{\lvert #1\rvert}
\begin{document}

\title{\bf\Large Clustered Variants of Haj\'os' Conjecture\footnote{This material is based upon work supported by the National Science Foundation under Grant No.\ DMS-1929851 and DMS-1954054.}}

\author{Chun-Hung Liu\footnote{Department of Mathematics, Texas A\&M University, College Station, Texas, USA, \texttt{chliu@math.tamu.edu}. Partially supported by NSF under Grant No.\ DMS-1929851 and DMS-1954054.} \quad David R. Wood\footnote{School of Mathematics, Monash University, Melbourne, Australia, \texttt{david.wood@monash.edu}. Supported by the Australian Research Council.}}

\maketitle

\begin{abstract}
Haj\'os conjectured that every graph containing no subdivision of the complete graph $K_{s+1}$ is properly $s$-colorable. 
This conjecture was disproved by Catlin. Indeed, the maximum chromatic number of such graphs is $\Omega(s^2/\log s)$. 
We prove that $O(s)$ colors are enough for a weakening of this conjecture that only requires every monochromatic component to have bounded size (so-called \emph{clustered} coloring).

Our approach leads to more results, many of which only require a much weaker assumption that forbids an `almost $(\leq 1)$-subdivision' (where at most one edge is subdivided more than once). 
This assumption is best possible, since no bound on the number of colors exists unless we allow at least one edge to be subdivided arbitrarily many times.
We prove the following (where $s \geq 2$):
\begin{enumerate}
	\item Graphs of bounded treewidth and with no almost $(\leq 1)$-subdivision of $K_{s+1}$ are $s$-choosable with bounded clustering.
	\item For every graph $H$, graphs with no $H$-minor and no almost $(\leq 1)$-subdivision of $K_{s+1}$ are $(s+1)$-colorable with bounded clustering.
	\item For every graph $H$ of maximum degree at most $d$, graphs with no $H$-subdivision and no almost $(\leq 1)$-subdivision of $K_{s+1}$ are $\max\{s+3d-5,2\}$-colorable with bounded clustering.
	\item For every graph $H$ of maximum degree $d$, graphs with no $K_{s,t}$ subgraph and no $H$-subdivision are $\max\{s+3d-4,2\}$-colorable with bounded clustering.
	\item Graphs with no $K_{s+1}$-subdivision are $(4s-5)$-colorable with bounded clustering. 
\end{enumerate}
The first result is tight and shows that the clustered analogue of Haj\'{o}s' conjecture is true for graphs of bounded treewidth. 
The second result implies an upper bound for the clustered version of Hadwiger's conjecture that is only one color away from the known lower bound, and shows that the number of colors is independent of the forbidden minor.  
The final result is the first $O(s)$ bound on the clustered chromatic number of graphs with no $K_{s+1}$-subdivision.
\end{abstract}

\newpage
\section{Introduction}

In the 1940s, Haj\'{o}s conjectured that every graph containing no subdivision of the complete graph $K_{s+1}$ is $s$-colorable; see \citep{Thomassen05,SeymourHC,Mohar05}. \citet{Dirac52} proved the conjecture for $s\leq 3$. 
It is open for $s \in \{4,5\}$, which would imply the Four Color Theorem. \citet{Catlin79}  presented counterexamples for all $s\geq 6$, and \citet{EF81} proved that the conjecture is false for almost all graphs. Indeed, there are graphs with no $K_{s+1}$-subdivision and with chromatic number $\Omega(s^2/\log s)$.  The best upper bound on the number of colors is $O(s^2)$, independently due to \citet{BT98} and \citet{KS96}; see \citep{FLS13} for a related result. See \citep{Mohar05,Thomassen05} for more explicit counterexamples and further discussion of connections to other areas of graph theory. 

The purpose of this paper is to prove several positive results in the direction of weakenings of {H}aj\'os' conjecture. 
Define a \emph{coloring} of a graph $G$ to simply be a function that assigns one color to each vertex of $G$. 
For a coloring $c$ of a graph $G$, a \emph{monochromatic $c$-component} of $G$ is a connected component of a subgraph of $G$ induced by all the vertices assigned the same color by $c$. When $c$ is clear, we simply write  \emph{monochromatic component}. 
A coloring has \emph{clustering} $\eta$ if every monochromatic component has at most $\eta$ vertices. Our focus is on minimizing the number of colors, with small clustering as a secondary goal. The \emph{clustered chromatic number} of a graph class $\mathcal{F}$ is the minimum integer $k$ for which there exists an integer $c$ such that every graph in $\mathcal{F}$ has a $k$-coloring with clustering $c$. There have been several recent papers on this topic \citep{NSSW19,vdHW18,KO19,ADOV03,CE19,Kawa08,KM07,LMST08,EJ14,HST03,EO16,DN17,LO17,HW19,MRW17}; see \citep{WoodSurvey} for a survey.

Most of our results actually hold (in some sense) for more general classes of graphs than those with no $K_{s+1}$-subdivision, as we now explain. 
Say a graph $H'$ is an {\it almost $(\leq 1)$-subdivision} of a graph $H$ if $H'$ can be obtained from $H$ by subdividing edges, where at most one edge is subdivided more than once. Most of our results say that all graphs containing no almost $(\leq 1)$-subdivision of $K_{s+1}$, plus some other properties, are $s$-colorable with bounded clustering. 

The following is our first main result. It provides a Haj\'os-type result for clustered coloring of graphs with bounded treewidth. 

\begin{theorem}
\label{TreewidthSubdivisionChoose}
For all $s,w \in \mathbb{N}$, there exists $\eta \in {\mathbb N}$ such that every graph with treewidth at most $w$ and with no almost $(\leq 1)$-subdivision of $K_{s+1}$ is $s$-choosable with clustering $\eta$. 
\end{theorem}

The notion of $s$-choosable with bounded clustering is defined in \cref{sec:listcoloring}. Note that every graph that is $s$-choosable with bounded clustering is also $s$-colorable with bounded clustering. This shows that the number of colors in \cref{TreewidthSubdivisionChoose} is best possible in the following strong sense: for all $s\in\mathbb{N}$ and $\eta \in \mathbb{N}$ there is a graph $G$ with treewidth at most $s-1$ (and thus with no subdivision of $K_{s+1}$), such that every  $(s-1)$-coloring of $G$ has a monochromatic component with at least $\eta$ vertices; see \citep{WoodSurvey}. In particular, at least $s$ colors are required even for this weakening of Haj\'{o}s' conjecture.

The assumption of bounded treewidth in  \cref{TreewidthSubdivisionChoose} is equivalent to saying that the graph excludes a planar graph as a minor by Robertson and Seymour's Grid Minor Theorem \citep{RS-V}. What if we exclude a general graph as a minor? Our next result answers this question (with one more color). 

\begin{theorem}
\label{MinorSubdivisionColor}
For every $s\in\mathbb{N}$ and every graph $H$, there exists $\eta \in {\mathbb N}$ such that every graph containing no $H$-minor and containing no almost $(\leq 1)$-subdivision of $K_{s+1}$ is  $(s+1)$-colorable with clustering $\eta$. 
\end{theorem}

\cref{MinorSubdivisionColor} (with $H=K_{s+1}$) has the following interesting corollary for graphs excluding a minor. 

\begin{corollary}
\label{MinorColor}
For every $s\in\mathbb{N}$ there exists $\eta \in {\mathbb N}$ such that every graph containing no $K_{s+1}$-minor is  $(s+1)$-colorable with  clustering $\eta$. 
\end{corollary}

\citet{KM07} first proved that graphs containing no $K_{s+1}$-minor are $O(s)$-colorable with bounded clustering. The bound on the number of colors has since been steadily improved \citep{Wood10,EKKOS15,LO17,Norin15,vdHW18}. 
Prior to the present work, the best bound was $s+2$, which followed from a general result by the authors \citep{LW2}. \cref{MinorColor} improves this bound to $s+1$, although it should be noted that results from \citep{LW2} are essential for the proof of \cref{MinorSubdivisionColor,MinorColor}. \citet{DN17} have announced that a forthcoming paper will prove that $s$ colors suffice (which is the clustered analogue of Hadwiger's Conjecture, and would be best possible). Their result is incomparable with \cref{MinorSubdivisionColor} and the aforementioned general result in \citep{LW2}.

Our next result relaxes the assumption that the graph contains no $H$-minor, and instead assumes that it contains no $H$-subdivision. The price paid is an increase in the number of colors, depending only on the maximum degree of $H$.

\begin{theorem}
\label{SubdivisionSubdivisionColor}
For every $s\in\mathbb{N}$ and every graph $H$ with maximum degree $d\in\mathbb{N}$, there exists $\eta\in\mathbb{N}$ such that every graph with no $H$-subdivision and no almost $(\leq 1)$-subdivision of $K_{s+1}$ is $\max\{s+3d-5,2\}$-colorable with clustering $\eta$. 
\end{theorem}

The next theorem relaxes the assumption of no almost $(\leq 1)$-subdivision of $K_{s+1}$, and instead assumes the graph contains no $K_{s,t}$-subgraph. Interestingly the number of colors does not depend on $t$. 
Note that $K_{s,t}$ contains a $K_{s+1}$-subdivision where every edge is subdivided at most once, when $t$ is sufficiently large.

\begin{theorem}
\label{SubdivisionSubgraphColor}
For $s,t,d \in \mathbb{N}$ and every graph $H$ of maximum degree $d$, there exists $\eta \in {\mathbb N}$ such that every graph with no $K_{s,t}$-subgraph and no $H$-subdivision is $\max\{s+3d-4,2\}$-colorable with clustering $\eta$. 
\end{theorem}

We remark that all of the above theorems forbid $(\leq 1)$-subdivisions of $K_{s+1}$ or subdivisions of $H$.
That is, we forbid a subdivision of a graph where some edge is allowed to be subdivided arbitrarily many times.
This condition is required since there are graphs of arbitrarily high girth and arbitrarily high chromatic number
\citep{Erdos59}, which therefore require arbitrarily many colors for any fixed clustering value; this shows that excluding finitely many graphs as subgraphs cannot ensure any upper bound on the number of colors.

Our final theorem simply excludes a $K_{s+1}$-subdivision. This is the first $O(s)$ bound on the clustered chromatic number of the class of graphs excluding a $K_{s+1}$-subdivision.

\begin{theorem}
\label{SubdivisionColor}
For each $s \in \mathbb{N}$, there exists $\eta \in {\mathbb N}$ such that every graph containing no $K_{s+1}$-subdivision is $\max\{4s-5,1\}$-colorable with clustering $\eta$. 
\end{theorem}

We now compare the above theorems with Haj\'os' conjecture. 
First note that  Theorems~\ref{TreewidthSubdivisionChoose}--\ref{SubdivisionSubdivisionColor} are stronger than Haj\'os' conjecture in the sense that they only exclude an almost $(\leq 1)$-subdivision of $K_{s+1}$, whereas  Haj\'os' conjecture excludes all subdivisions of $K_{s+1}$. 
Moreover, \cref{TreewidthSubdivisionChoose} also holds in the stronger setting of choosability. On the other hand, Theorems~\ref{TreewidthSubdivisionChoose}--\ref{SubdivisionColor} are weaker than Haj\'os' conjecture in the sense that they have bounded clustering rather than a proper coloring. However, such a weakening is unavoidable since Haj\'os' conjecture is false. Indeed, the proof of the theorem of \citet{EF81} mentioned above shows that, for a suitable constant $c$, almost every graph on $cs^2$ vertices contains no subdivision of $K_{s+1}$ and has chromatic number $\Omega(s^2/\log s)$. Trivially, such a graph has treewidth at most $cs^2$ and contains no $K_{cs^2}$-minor. Thus the clustering function in all of the above theorems is at least $\Omega(s/\log s)$.

The paper is organized as follows. \cref{Preliminaries} introduces preliminary definitions and results from our companion papers \citep{LW1,LW2} that are used in the present paper. 
\cref{ExcludingTopologicalMinors} introduces a structure theorem of the first author and Thomas~\cite{LT19} for graphs excluding a fixed subdivision, and uses it to prove \cref{SubdivisionSubgraphColor}. Building on this work, \cref{OneStepFurther} proves the remaining theorems mentioned above.

%%%%%%%%%%%%%%%%%
\section{Preliminaries}
\label{Preliminaries}

We use the following notation. 
Let 
$\mathbb{N}_0:=\{0,1,2,\dots\}$ and
$\mathbb{N}:=\{1,2,\dots\}$. 
For $m,n\in\mathbb{N}_0$, let $[m,n]:=\{m,m+1,\dots,n\}$ and $[n]:=[1,n]$. 

Let $G$ be a graph (allowing loops and parallel edges). For $v\in V(G)$, let $N_G(v):=\{w\in V(G): vw\in E(G)\}$ be the neighborhood of $v$, and let  $N_G[v] :=N_G(v) \cup\{v\}$. For $X\subseteq V(G)$, let $N_G(X):= \bigcup_{v\in X} (N_G(v) - X)$ and $N_G[X]:= N_G(X)\cup X$. Denote the subgraph of $G$ induced by $X$ by $G[X]$.

For a graph $G$, a subset $X$ of $V(G)$, and an integer $s\geq 1$, let
\begin{align*}
\ngs[G]{X} & := \{v \in V(G)- X:  \lvert N_G(v) \cap X \rvert \geq s\} \text{ and} \\
\nls[G]{X} & := \{v \in V(G)- X:  1 \leq \lvert N_G(v) \cap X \rvert < s\}.
\end{align*}
When the graph $G$ is clear from the context we write $\ngs{X}$ instead of $\ngs[G]{X}$, and similarly for 
$\nls{X}$.

\begin{lemma}[{\protect\cite[Lemma~12]{LW1}}]
\label{BoundedGrowth}
For all $s,t\in\mathbb{N}$, there exists a function $f_{s,t}:\mathbb{N}_0\rightarrow\mathbb{N}_0$ such that for every graph $G$ with no $K_{s,t}$ subgraph, if $X\subseteq V(G)$ then $|\ngs{X}| \leq f_{s,t}(\lvert X \rvert)$.
\end{lemma}

\cref{BoundedGrowth} is sufficient to prove the theorems in this paper. But when $G$ excludes a fixed minor or subdivision of a fixed graph, the function $f_{s,t}$ in \cref{BoundedGrowth} can be made linear; see \citep{LW1}. This improves the clustering function in all our results, although to simplify the presentation, we choose not to explicitly evaluate our clustering functions.

A \emph{tree-decomposition} of a graph $G$ is a pair $(T,\mathcal{X}=(X_x:x\in V(T)))$, where $T$ is a tree and for each node $x\in V(T)$,  $X_x$ is a subset of $V(G)$ called a \emph{bag}, such that for each vertex $v\in V(G)$, the set $\{x\in V(T):v\in X_x\}$ induces a non-empty (connected) subtree of $T$, and for each edge $vw\in E(G)$ there is a node $x\in V(T)$ such that $\{v,w\} \subseteq X_x$. 
The \emph{width} of a tree-decomposition $(T,\mathcal{X})$ is $\max\{|X_x|-1: x\in V(T)\}$. The \emph{treewidth} of a graph $G$ is the minimum width of a tree-decomposition of $G$. 

Let $H$ be a graph. An {\it $H$-minor} of a graph $G$ is a map $\alpha$ with domain $V(H) \cup E(H)$ such that:
\begin{itemize}
	\item For every $h \in V(H)$, $\alpha(h)$ is a nonempty connected subgraph of~$G$. 
	\item If $h_1$ and $h_2$ are different vertices of~$H$, then $\alpha(h_1)$ and $\alpha(h_2)$ are disjoint.
	\item For each edge $e$ of~$H$ with endpoints $h_1,h_2$, $\alpha(e)$ is an edge of~$G$ with one end in $\alpha(h_1)$ and one end in $\alpha(h_2)$; furthermore, if $h_1=h_2$, then $\alpha(e) \in E(G)-E(\alpha(h_1))$.
        \item If $e_1, e_2$ are two different edges of~$H$, then $\alpha(e_1) \neq \alpha(e_2)$.
\end{itemize}

\subsection{List Coloring} \label{sec:listcoloring}

For our purposes, a \emph{color} is an element of $\mathbb{Z}$. 
Let $G$ be a graph.
A \emph{list-assignment} of $G$ is a function $L$ with domain containing $V(G)$, such that $L(v)$ is a non-empty set of colors for each vertex $v\in V(G)$. 
For a list-assignment $L$ of $V(G)$, an {\it $L$-coloring} of $G$ is a coloring $c$ of $G$ such that $c(v) \in L(v)$ for every $v \in V(G)$. 
An $L$-coloring has \emph{clustering} $\eta$ if every monochromatic component has at most $\eta$ vertices. 
A list-assignment $L$ of  $G$ is an \emph{$\ell$-list-assignment} if $|L(v)|\geq\ell$ for every vertex $v\in V(G)$. 
A graph $G$ is \emph{$\ell$-choosable with clustering $\eta$} if $G$ is $L$-colorable with clustering $\eta$ for every $\ell$-list-assignment $L$ of $G$. 

For a graph $G$, a subset $Y_1\subseteq V(G)$, and  $s,r\in\mathbb{N}$, a list-assignment $L$ of $G$ is an 
\emph{$(s,r,Y_1)$-list-assignment} if:
	\begin{itemize}
		\item[(L1)] $ \lvert L(v) \rvert \in [ s+r]$ for every $v \in V(G)$.
		\item[(L2)] $Y_1 = \{v \in V(G): \lvert L(v) \rvert=1\}$. 
		\item[(L3)] For every $y \in \nls{Y_1}$, 
			$$\lvert L(y) \rvert = s+r-\lvert N_G(y) \cap Y_1 \rvert \geq r+1$$ and $L(y) \cap L(u)=\emptyset$ for every $u \in N_G(y) \cap Y_1$.
		\item[(L4)] For every $v \in V(G)-N_G[Y_1]$, we have $\lvert L(v) \rvert =s+r$.
		\item[(L5)] For every $v \in V(G)-Y_1$, we have $\lvert L(v) \rvert \geq r+1$.
	\end{itemize}
We say that an $(s,r,Y_1)$-list-assignment $L$ of $G$ is \emph{restricted} if:
\begin{itemize}
\item[(L1')] $ L(v) \subseteq [s+r]$ for every $v \in V(G)$.
\end{itemize}
Note that a restricted $(s,2,Y_1)$-list-assignment is called a $(s,Y_1,0,0)$-list-assignment in our companion paper \citep{LW2}. 

For  a list-assignment $L$ of a graph $G$ with $Y_1=\{v \in V(G): \lvert L(v) \rvert=1\}$,  
for $\eta\in\mathbb{N}$ and a nondecreasing function $g:\mathbb{N}\to\mathbb{N}$, 
an $L$-coloring $c$ of $G$ is {\it $(\eta,g)$-bounded} if:
	\begin{itemize}
		\item the union of the monochromatic components intersecting $Y_1$ contains at most $\lvert Y_1 \rvert^2 g(\lvert Y_1 \rvert)$ vertices, and 
		\item every monochromatic component contains at most $\eta^2 g(\eta)$ vertices.
	\end{itemize}

\subsection{Companion Results}

Our companion paper proves the following results for graphs with no $K_{s,t}$ subgraph. The first assumes  bounded treewidth, the second assumes an excluded minor.

\begin{theorem}[{\protect\citep[Theorem~17]{LW2}}]
\label{bdd tw Kst}
For all $s,t,w\in\mathbb{N}$, there exist $\eta\in\mathbb{N}$ and a nondecreasing function $g$ such that if $G$ is a graph of treewidth at most $w$ and with no $K_{s,t}$ subgraph, $Y_1$ is a subset of $V(G)$ with $\lvert Y_1 \rvert \leq \eta$, and $L$ is an $(s,1,Y_1)$-list-assignment of $G$, then there exists an $(\eta,g)$-bounded $L$-coloring of $G$.
\end{theorem}

\begin{theorem}[{\protect\citep[Theorem~24]{LW2}}]
\label{odd minor free}
For all $s,t\in\mathbb{N}$ and for every graph $H$, there exist $\eta\in\mathbb{N}$ and a nondecreasing function $g$ such that  if $G$ is a graph with no $K_{s,t}$ subgraph and no $H$-minor, $Y_1$ is a subset of $V(G)$ with $\lvert Y_1 \rvert \leq \eta$, and $L$ is a restricted $(s,2,Y_1)$-list-assignment of $G$, then there exists an $(\eta,g)$-bounded $L$-coloring.
\end{theorem}

\subsection{Progress}

The concept of ``progress'' from the proofs of the above two theorems are re-used in the present paper. 
Let $s,r \in \mathbb{N}$ and $L$ be an $(s,r,Y_1)$-list-assignment of a graph $G$. 
For $W \subseteq V(G)$, a {\it $W$-progress} of $L$ is a list-assignment $L'$ of $G$ defined as follows: 
\begin{itemize}
	\item Let $Y_1':=Y_1 \cup W$. 
	\item For every $y \in Y_1$, let $L'(y):=L(y)$.
	\item For every $y \in Y_1'-Y_1$, let $L'(y)$ be a 1-element subset of $L(y)$ (which exists by (L2)--(L5)). 
	\item For each $v \in \nls{Y_1'}$, let $L'(v)$  be a subset of $L(v)-\{L'(w): w \in N_G(v) \cap (W-Y_1)\}$ of size $\lvert L(v) \rvert - \lvert N_G(v) \cap (W-Y_1) \rvert$.
	\item For every $v \in V(G)-(Y_1' \cup \nls{Y_1'})$, let $L'(v):=L(v)$.
\end{itemize}

\begin{lemma}[{\protect\citep[Lemma~12~with~$F=\emptyset$]{LW2}}] 
 \label{progress basic}
  Let $s,r \in \mathbb{N}$ and $L$ be an $(s,r,Y_1)$-list-assignment of a graph $G$. 
Let $W \subseteq V(G)$. Then every $W$-progress $L'$ of $L$  is an $(s,r,Y_1 \cup W)$-list-assignment of $G$, and  $L'(v) \subseteq L(v)$ for every $v \in V(G)$.
\end{lemma}

\begin{lemma}[{\protect\citep[Lemma~13 with $F=\emptyset$]{LW2}}] 
\label{enlarge Y}
For all $s,t,k\in\mathbb{N}$, there exist a number $\eta>k$ and a nondecreasing function $g$ with domain ${\mathbb N_0}$ and with $g(0) \geq \eta$ such that if $G$ is a graph with no $K_{s,t}$ subgraph, $r \in \mathbb{N}$, $Y_1$ is a subset of $V(G)$ with $\lvert Y_1 \rvert \leq \eta$, and $L$ is an $(s,r,Y_1)$-list-assignment of $G$, then at least one of the following holds:
	\begin{enumerate}
		\item There exists an $(\eta,g)$-bounded $L$-coloring of $G$.
		\item $\lvert Y_1 \rvert >k$.
		\item For every color $\ell$, there exist a subset $Y_1'$ of $V(G)$ with $\eta \geq \lvert Y_1' \rvert > \lvert Y_1 \rvert$ and an $(s,r,Y_1')$-list-assignment $L'$ of $G$ with $L'(v) \subseteq L(v)$ for every $v \in V(G)$, such that:
			\begin{enumerate}
				\item there does not exist an $(\eta,g)$-bounded $L'$-coloring $c'$ of $G$,
				\item for every $L'$-coloring of $G$, every monochromatic component intersecting $Y_1$ is contained in $G[Y_1']$, and
				\item for every $y \in Y_1'$ with $\ell \in L'(y)$, we have $\{v \in N_G(y)-Y_1': \ell \in L'(v)\}=\emptyset$.
			\end{enumerate}
		\item $Y_1 \neq \emptyset$, $N_G(Y_1)=\emptyset$, and there does not exist an $(\eta,g)$-bounded $L|_{G-Y_1}$-coloring of $G-Y_1$. 
	\end{enumerate}
\end{lemma}

\subsection{Separations and Tangles}

A {\it separation} of a graph $G$ is an ordered pair $(A,B)$ of edge-disjoint subgraphs of $G$ with $A \cup B=G$. 
The {\it order} of $(A,B)$ is $\lvert V(A \cap B) \rvert$.
A {\it tangle} $\T$ in a graph $G$ of order $\theta$ is a set of separations of $G$ of order less than $\theta$ such that:
	\begin{itemize}
		\item[(T1)] For every separation $(A,B)$ of $G$ of order less than $\theta$, either $(A,B) \in \T$ or $(B,A) \in \T$.
		\item[(T2)] If $(A_i,B_i) \in \T$ for $i \in [3]$, then $A_1 \cup A_2 \cup A_3 \neq G$.
		\item[(T3)] If $(A,B) \in \T$, then $V(A) \neq V(G)$.
	\end{itemize}

\begin{lemma}[{\protect\citep[Lemma~16 with $F=\emptyset$]{LW2}}]
\label{coloring or tangle}
For all $s,t,\theta,\eta,r\in\mathbb{N}$ with $\eta \geq 9\theta+1$, for every nondecreasing function $g$ with domain ${\mathbb N_0}$, if $G$ is a graph with no $K_{s,t}$ subgraph, $Y_1$ is a subset of $V(G)$ with $9\theta+1 \leq \lvert Y_1 \rvert \leq \eta$, and $L$ is an $(s,r,Y_1)$-list-assignment of $G$, then at least one of the following holds:
	\begin{enumerate}
		\item There exists an $(\eta,g)$-bounded $L$-coloring of $G$.
		\item There exist an induced subgraph $G'$ of $G$ with $\lvert V(G') \rvert < \lvert V(G) \rvert$, a subset $Y_1'$ of $V(G')$ with $\lvert Y_1' \rvert \leq \eta$ and an $(s,r,Y_1')$-list-assignment $L'$ of $G'$ such that:
		\begin{enumerate}
			\item $L'(v) \subseteq L(v)$ for every $v \in V(G')$.
			\item There does not exist an $(\eta,g)$-bounded $L'$-coloring of $G'$.
		\end{enumerate} 
		\item $\T:=\{(A,B): \lvert V(A \cap B) \rvert < \theta, \lvert V(A) \cap Y_1 \rvert \leq 3\theta\}$ is a tangle of order $\theta$ in $G$.
	\end{enumerate}
\end{lemma}

A tangle~$\T$ in $G$ \emph{controls} an $H$-minor $\alpha$ if there does not exist $(A,B) \in \T$ of order less than $\abs{V(H)}$ such that $V(\alpha(h)) \subseteq V(A)$ for some $h\in V(H)$.

\begin{lemma}[{\protect\citep[Lemma~23 with $\ell=r=0$]{LW2}}] 
\label{not controlling minor}
	For all $s,t,t' \in\mathbb{N}$, there exist $\theta^* \in \mathbb{N}$ and nondecreasing functions $g^*,\eta^*$ with domain ${\mathbb N_0}$ such that if $G$ is a graph with no $K_{s,t}$ subgraph, $\theta \in \mathbb{N}$ with $\theta \geq \theta^*$, 
	$\eta \in \mathbb{N}$ with $\eta \geq \eta^*(\theta)$, $Y_1 \subseteq V(G)$ with $3\theta < \lvert Y_1 \rvert \leq \eta$, 
	$L$ is a restricted $(s,2,Y_1)$-list-assignment of $G$, 
	$g$ is a nondecreasing function with domain ${\mathbb N_0}$ with $g \geq g^*$, and 
	$\T := \{(A,B): \lvert V(A \cap B) \rvert < \theta, \lvert V(A) \cap Y_1 \rvert \leq 3\theta\}$ 
	is a tangle in $G$ of order $\theta$ that does not control a $K_{t'}$-minor, then either:
	\begin{enumerate}
		\item there exists an $(\eta,g)$-bounded $L$-coloring of $G$, or 
		\item there exist $(A^*,B^*) \in \T$, a set $Y_{A^*}$ with $\lvert Y_{A^*} \rvert \leq \eta^*(\theta)$ and $Y_1 \cap V(A^*) \subseteq Y_{A^*} \subseteq V(A^*)$, and a restricted $(s,2,Y_{A^*})$-list-assignment $L_{A^*}$ of $G[V(A^*)]$ such that there exists no $(\eta,g)$-bounded $L_{A^*}$-coloring of $G[V(A^*)]$.
	\end{enumerate}
\end{lemma}

%%%%%%%%%%%%%%%%%%%%%%
\section{Excluding Subdivisions}
\label{ExcludingTopologicalMinors}

The following theorem is a special case of a theorem by the first author and Thomas~\cite{LT19}.

\begin{theorem}[{\cite[Theorem~6.8]{LT19}}] \label{StructureExcludingSubdivision}
For any integers $d,h$ and graph $H$ on $h$ vertices with maximum degree at most $d$, there exist integers $\theta,\xi$ such that if $G$ is a graph containing no $H$-subdivision, and if $\T$ is a tangle in $G$ of order at least $\theta$ controlling a $K_{\lfloor \frac{3}{2}dh \rfloor}$-minor, then there exists $Z \subseteq V(G)$ with $\lvert Z \rvert \leq \xi$ such that for every vertex $v \in V(G)-Z$, there exists $(A,B) \in \T-Z$ of order less than $d$ such that $v \in V(A)-V(B)$.
\end{theorem}

The next two lemmas imply \cref{SubdivisionSubgraphColor}, since if $s,d \in {\mathbb N}$ and $3d+s<7$, then $d=1$.

\begin{lemma} \label{easyforbidmaxdeg1subdiv}
If $H$ is a graph of maximum degree at most 1, then every graph with no $H$-subdivision is 2-colorable with clustering  $\max\{2\lvert V(H) \rvert-2,1\}$.
\end{lemma}

\begin{proof}
Since $H$ is of maximum degree at most one, $G$ has no $H$-subdivision implies that $G$ does not contain a matching of size $\lvert V(H) \rvert$, and hence $G$ contains a vertex-cover $S$ of size at most $2\lvert V(H) \rvert-2$.
By coloring every vertex in $S$ with 1 and coloring every vertex in $V(G)-S$ with 2, we obtain a 2-coloring of $G$ with clustering $\max\{\lvert S \rvert,1\} \leq \max\{2\lvert V(H) \rvert-2,1\}$.
\end{proof}

\begin{lemma} 
\label{no H subdivision Kst subgraph}
For any $s,t,d \in \mathbb{N}$ and graph $H$ of maximum degree $d$ with $3d+s \geq 7$, there exist $\eta\in\mathbb{N}$ and a nondecreasing function $g$ such that if $G$ is a graph with no $K_{s,t}$ subgraph and no $H$-subdivision, $Y_1 \subseteq V(G)$ with $\lvert Y_1 \rvert \leq \eta$ and $L$ is a restricted $(s',2,Y_1)$-list-assignment of $G$, then there exists an $(\eta,g)$-bounded $L$-coloring, where $s'=3d+s-6$. 
\end{lemma}

\begin{proof}
Define the following:
	\begin{itemize}
		\item Let $f$ be the function $f_{s,t}$ mentioned in \cref{BoundedGrowth}.
		\item Let $\theta_0$ be the number $\theta^*$ and $g_0,\eta_0$ be the functions $g^*,\eta^*$, respectively, mentioned in \cref{not controlling minor} by taking $s=s', t=t$ and $t'=\lfloor \frac{3}{2}d\lvert V(H) \rvert \rfloor$.
		\item Let $\theta_1$ and $\xi$ be the numbers $\theta$ and $\xi$ mentioned in \cref{StructureExcludingSubdivision}, respectively, by taking $d=d$, $h=\lvert V(H) \rvert$ and $H=H$.
		\item Let $a_0:=f(\xi)d^2 + \xi+1$, and let $a_i:=da_{i-1}+1$ for $i \in\mathbb{N}$. 
		\item Let $\theta := \theta_0+\theta_1 + (d-1)a_{(d-1)a_0}$.
		\item Let $\eta_1$ be the number $\eta$ and let $g_1$ be the function $g$ mentioned in \cref{enlarge Y} by taking $s=s'$, $t=t$ and $k=9\theta$.
		Note that $g(0) \geq \eta_1 > 9\theta$ by \cref{enlarge Y}.
		\item Let $\eta := \eta_0(\theta) + \eta_1 + (d-1)a_{(d-1)a_0}$.
		\item Let $g:\mathbb{N}\to\mathbb{N}$ be the function defined by $g(0):=g_0(0)+g_1(0)$ and $g(x+1) := g_0(x+1) + g_1(x+1)+ \sum_{i=0}^x i^2g(i)$ for $x\in\mathbb{N}$.
	\end{itemize}

Let $G$ be a graph with no $K_{s,t}$ subgraph and with no subdivision of $H$, let $Y_1 \subseteq V(G)$ with $\lvert Y_1 \rvert \leq \eta$, and let $L$ be a restricted $(s',2,Y_1)$-list-assignment of $G$.
Suppose to the contrary that there exists no $(\eta,g)$-bounded $L$-coloring of $G$.
We further assume that $\lvert V(G) \rvert$ is minimum, and subject to this, $\lvert Y_1 \rvert$ is maximum.

\begin{claim}
\label{no H subdivision Kst subgraph Claim1}
$Y_1 \neq \emptyset$ and $N_G(Y_1) \neq \emptyset$.
\end{claim}

\begin{proof} 
First suppose that $Y_1=\emptyset$.
Let $v$ be a vertex of $G$, and let $L'$ be a $\{v\}$-progress of $L$.
Let $Y_1'=\{v\}$.
By \cref{progress basic}, $L'$ is an $(s',2,Y_1')$-list-assignment of $G$.
Since $\lvert Y_1' \rvert \leq \eta$, the maximality of $\lvert Y_1 \rvert$ implies that there exists an $(\eta,g)$-bounded $L'$-coloring $c'$ of $G$.
But $c'$ is an $(\eta,g)$-bounded $L$-coloring $c$ of $G$, a contradiction.

So $Y_1 \neq \emptyset$.
Suppose that $N_G(Y_1)=\emptyset$.
Let $G':=G-Y_1$.
Then $L|_{G'}$ is an $(s',2,\emptyset)$-list-assignment of $G-Y_1$. 
By the minimality of $\lvert V(G) \rvert$, there exists an $(\eta,g)$-bounded $L|_{G'}$-coloring $c$ of $G'$.
Color each vertex $y$ in $Y_1$ with the unique element in $L(y)$. 
Since $\lvert Y_1 \rvert \leq \lvert Y_1 \rvert^2g(\lvert Y_1 \rvert)$, we obtain an $(\eta,g)$-bounded $L$-coloring of $G$, a contradiction.
\end{proof}

\begin{claim}
\label{no H subdivision Kst subgraph Claim2}
$\lvert Y_1 \rvert \geq 9\theta+1$.
\end{claim}

\begin{proof} 
Suppose $\lvert Y_1 \rvert \leq 9\theta$.
So $\lvert Y_1 \rvert <\eta_1$.
Since $G$ has no $K_{s,t}$ subgraph, $G$ has no $K_{s',t}$ subgraph.
Applying \cref{enlarge Y,no H subdivision Kst subgraph  Claim1}, either there exists an $(\eta_1,g_1)$-bounded $L$-coloring of $G$, or there exist $Y_1' \subseteq V(G)$ with $\eta_1 \geq \lvert Y_1' \rvert > \lvert Y_1 \rvert$ and an $(s',2,Y_1')$-list-assignment $L'$ of $G$ with $L'(v) \subseteq L(v)$ for every $v \in V(G)$ such that for every $L'$-coloring of $G$, every monochromatic component intersecting $Y_1$ is contained in $G[Y_1']$.
Since $\eta_1 \leq \eta$ and $g_1 \leq g$, every $(\eta_1,g_1)$-bounded $L$-coloring of $G$ is an $(\eta,g)$-bounded $L$-coloring of $G$, so the former does not hold.
Hence there exist $Y_1' \subseteq V(G)$ with $\eta_1 \geq \lvert Y_1' \rvert > \lvert Y_1 \rvert$ and a restricted $(s',2,Y_1')$-list-assignment $L'$ of $G$ with $L'(v) \subseteq L(v)$ for every $v \in V(G)$ such that for every $L'$-coloring of $G$, every monochromatic component intersecting $Y_1$ is contained in $G[Y_1']$.
Since $\lvert Y_1' \rvert \leq \eta_1 \leq \eta$, the maximality of $\lvert Y_1 \rvert$ implies that there exists an $(\eta,g)$-bounded $L'$-coloring $c'$ of $G$.
So every monochromatic component respect to $c'$ contains at most $\eta^2g(\eta)$ vertices.
Since $L'(v) \subseteq L(v)$ for every $v \in V(G)$, $c'$ is also an $L$-coloring of $G$.
Every monochromatic $c'$-component intersecting $Y_1$ is contained in $G[Y_1']$ and hence contains at most $\lvert Y_1' \rvert \leq \eta_1 \leq g_1(0) \leq g(0) \leq \lvert Y_1 \rvert^2g(\lvert Y_1 \rvert)$ vertices.
So $c'$ is an $(\eta,g)$-bounded $L$-coloring of $G$, a contradiction.
\end{proof}

Let $\T$ be the set of separations $(A,B)$ of $G$ such that $\lvert V(A \cap B) \rvert <\theta$ and $\lvert V(A) \cap Y_1 \rvert \leq 3\theta$. 

\begin{claim}
\label{no H subdivision Kst subgraph Claim3}
$\T$ is a tangle in $G$ of order $\theta$.
\end{claim}

\begin{proof} 
Suppose that $\T$ is not a tangle in $G$ of order $\theta$.
Note that $G$ has no $K_{s',t}$ subgraph and $L$ is an $(s',2,Y_1)$-list-assignment of $G$ with $\eta \geq \lvert Y_1 \rvert \geq 9\theta+1$ by  \cref{no H subdivision Kst subgraph Claim2}.
Applying \cref{coloring or tangle} by taking $s=s'$, $t=t$, $\theta=\theta$, $\eta=\eta$, $r=2$ and $g=g$, there exists an induced subgraph $G'$ of $G$ with $\lvert V(G') \rvert < \lvert V(G) \rvert$, a subset $Y_1' \subseteq V(G')$ with $\lvert Y_1' \rvert \leq \eta$, and an $(s',2,Y_1')$-list-assignment $L'$ of $G'$ with $L'(v) \subseteq L(v)$ for every $v \in V(G)$ such that there exists no $(\eta,g)$-bounded $L'$-coloring of $G'$. This contradicts the minimality of $\lvert V(G) \rvert$.
\end{proof}

\begin{claim}
\label{no H subdivision Kst subgraph Claim4}
$\T$ controls a $K_{\lfloor \frac{3}{2} d \lvert V(H) \rvert \rfloor}$-minor.
\end{claim}

\begin{proof}
Suppose to the contrary that $\T$ does not control a $K_{\lfloor \frac{3}{2} d \lvert V(H) \rvert \rfloor}$-minor.
Note that $\theta \geq \theta_0$, $\eta \geq \eta_0(\theta)$ and $g \geq g_0$.
Apply \cref{not controlling minor} with $s=s'$, $t=t$ and $t'=\lfloor \frac{3}{2} d \lvert V(H) \rvert \rfloor$. 
Since there does not exist an $(\eta,g)$-bounded $L$-coloring of $G$, we know there exist $(A^*,B^*) \in \T$, a set $Y_{A^*}$ with $\lvert Y_{A^*} \rvert \leq \eta_0(\theta) \leq \eta$ and $Y_1 \cap V(A^*) \subseteq Y_{A^*} \subseteq V(A^*)$, and a restricted $(s',2,Y_{A^*})$-list-assignment $L_{A^*}$ of $G[V(A^*)]$ such that there exists no $(\eta,g)$-bounded $L_{A^*}$-coloring of $G[V(A^*)]$.
But $\lvert V(A^*) \rvert < \lvert V(G) \rvert$ since $\lvert V(A^*) \cap Y_1 \rvert < 3\theta < \lvert Y_1 \rvert$.
This contradicts the minimality of $\lvert V(G) \rvert$.
\end{proof}

Since $G$ contains no subdivision of $H$, by \cref{StructureExcludingSubdivision} and \cref{no H subdivision Kst subgraph Claim4}, there exists $Z \subseteq V(G)$ with $\lvert Z \rvert \leq \xi$ such that for every $v \in V(G)-Z$, there exists $(A_v,B_v) \in \T-Z$ of order at most $d-1$ such that $v \in V(A_v)-V(B_v)$.

We may assume that for every $v \in V(G)-Z$,
	\begin{itemize}
		\item[(i)] $(A_v,B_v) \in \T-Z$ has order at most $d-1$ and $v \in V(A_v)-V(B_v)$,
		\item[(ii)] subject to (i), $A_v-V(A_v \cap B_v)$ is connected, 
		\item[(iii)] subject to (i) and (ii), every vertex in $V(A_v \cap B_v)$ is adjacent to some vertex in $V(A_v)-V(B_v)$, 
		\item[(iv)] subject to (i)--(iii), $V(A_v)$ is maximal,
		\item[(v)] subject to (i)--(iv), $\lvert V(A_v \cap B_v) \rvert$ is minimal, and 
		\item[(vi)] subject to (i)--(v), $A_v$ is maximal. 
	\end{itemize}
Note that for every $v \in V(G)-Z$, $A_v$ is connected and for every two vertices $x,y \in V(A_v)$, there exists a path in $A_v$ from $x$ to $y$ internally disjoint from $V(A_v \cap B_v)$ since $A_v - V(A_v \cap B_v)$ is connected and every vertex in $V(A_v \cap B_v)$ is adjacent to some vertex in $V(A_v)-V(B_v)$

For any subset $\C\subseteq \T-Z$, let $(A_\C,B_\C)$ be the separation $(\bigcup_{(A,B)\in\C}A, \bigcap_{(A,B)\in\C}B)$. 
Note that $V(A_\C \cap B_\C) \subseteq \bigcup_{(A,B) \in \C}V(A \cap B)$, so $\lvert V(A_\C \cap B_\C) \rvert \leq \lvert \C \rvert(d-1)$.

\begin{claim}
\label{no H subdivision Kst subgraph Claim5}
Let $\C=\{(A_w,B_w): w \in W\}$ for some $W \subseteq V(G)-Z$. 
If $x$ is a vertex in $V(A_\C \cap B_\C)$, then $V(A_x \cap B_x)-V(B_w) \neq \emptyset$ for some $w \in V(G)-Z$ with $(A_w,B_w) \in \C$.
\end{claim}

\begin{proof} 
Since $x \in V(A_\C \cap B_\C)$, there exists $w \in W \subseteq V(G)-Z$ such that $(A_w,B_w) \in \C$ and $x \in V(A_w \cap B_w)$.
Suppose to the contrary that $V(A_x \cap B_x) \subseteq V(B_w)$.

First suppose that there exists $v \in V(A_w)-(V(B_w) \cup V(A_x))$.
Since $A_w-V(B_w)$ is connected by (ii) and every vertex in $V(A_w \cap B_w)$ is adjacent to a vertex in $V(A_w)-V(B_w)$ by (iii), there exists a path $P$ in $G[(V(A_w)-V(B_w)) \cup \{x\}]$ from $x$ to $v$.
Since $x \in V(A_x)-V(B_x)$ and $v \in V(G)-(Z \cup V(A_x))$, $P-x$ intersects $V(A_x \cap B_x) \subseteq V(B_w)$.
But $V(P-x) \subseteq V(A_w)-V(B_w)$, a contradiction.
So $V(A_w)-V(B_w) \subseteq V(A_x)$.

Suppose that there exists a vertex $u \in V(A_w \cap B_w)-V(A_x)$.
Since $u \in V(A_w \cap B_w)$, there exists $u' \in N_G(u) \cap V(A_w)-V(B_w)$ by (iii).
So $u' \in N_G(u) \cap V(A_w)-V(B_w) \subseteq N_G(u) \cap V(A_x)$.
Since $u \not \in V(A_x)$, $u' \in V(A_x \cap B_x) \cap V(A_w)-V(B_w)$, contradicting the assumption $V(A_x \cap B_x) \subseteq V(B_w)$.
Hence $V(A_w \cap B_w) \subseteq V(A_x)$.

Therefore, $V(A_w) \subseteq V(A_x)$.
By (v), every vertex in $V(A_x \cap B_x)$ is adjacent to some vertex in $V(B_x)-V(A_x)$.
So if $V(A_x)=V(A_w)$, then $V(A_x \cap B_x) \subseteq V(A_w \cap B_w)$, and since $(A_w,B_w)$ satisfies (v), $V(B_w)=V(B_x)$.
Hence if $V(A_x)=V(A_w)$, then $(A_x,B_x)=(A_w,B_w)$ by (vi).
Since $x \in V(A_x)-V(B_x)$ and $x \in V(A_w \cap B_w)$, $(A_x,B_x) \neq (A_w,B_w)$.
So $V(A_w) \subset V(A_x)$.
Since $(A_w,B_w)$ satisfies (iv), $w \in V(B_x)$.
Since $V(A_x \cap B_x) \subseteq V(B_w)$ and $w \not \in V(B_w)$, $w \in V(B_x)-V(A_x)$.
So $V(A_w) \not \subseteq V(A_x)$, a contradiction.
\end{proof}

Let $Z':=\{v \in V(G)-(Y_1 \cup Z): \lvert N_G(v) \cap Z \rvert \geq s\}$.
Note that $\lvert Z' \rvert \leq f(\lvert Z \rvert) \leq f(\xi)$ by \cref{BoundedGrowth}.

We say that a triple $(\C,S,T)$ is {\it useful} if the following hold: 
\begin{itemize}
	\item[(U1)] There exists $W \subseteq V(G)-Z$ such that $\C=\{(A_v,B_v): v \in W\}$.
	\item[(U2)] $N_G[N_G[Z']] \cap V(B_\C)=\emptyset$. 
	\item[(U3)] $S$ is a subset of $N_G[V(A_\C \cap B_\C)] \cap V(A_\C)$ and $T$ is a subset of $Y_1 \cap V(A_\C)-V(B_\C)$ such that there exists a bijection $\iota$ from a subset of $Y_1 \cap V(A_\C)$ to $S$ such that:
	\begin{itemize}
		\item $\lvert S \rvert+ \lvert T \rvert + \lvert Z \rvert +1 \leq \lvert Y_1 \cap V(A_\C)-V(B_\C) \rvert+\lvert \{y \in Y_1 \cap S \cap V(A_\C \cap B_\C): \iota(y)=y\} \rvert$, and 
 		\item for every vertex $y$ in the domain of $\iota$,
		\begin{itemize}
			\item if $y \in V(A_\C)-V(B_\C)$ and there exists a vertex $v \in N_G(y) \cap V(A_\C \cap B_\C)-S$, then $\iota(y) \in N_G(y) \cap V(A_\C \cap B_\C)$, and 
			\item if $y \in V(A_\C \cap B_\C)$, then $\iota(y)=y$.
		\end{itemize}
	\end{itemize}
	\item[(U4)] $T$ is disjoint from $Z'$ and the domain of $\iota$.
\end{itemize}

\begin{claim}
\label{no H subdivision Kst subgraph Claim6}
There exists a collection $\C$ of members of $\T-Z$ with $\lvert \C \rvert \leq \lvert Z' \rvert d^2+\lvert Z \rvert +1$ such that $(\C,\emptyset,\emptyset)$ is useful. 
\end{claim}

\begin{proof} 
For every $u \in V(G)-Z$, let $\C_u := \{(A_u,B_u), (A_v,B_v): v \in N_G(u) \cap V(B_u)\}$.
Note that $\lvert N_G(u) \cap V(B_u) \rvert \leq \lvert V(A_u \cap B_u) \rvert \leq d-1$ since $u \in V(A_u)-V(B_u)$.
So $\lvert \C_u \rvert \leq d$.
Note that $N_G[\{u\}] \cap V(B_{\C_u}) = \emptyset$.

For every $u \in V(G)-Z$, let $\C_u' := \C_u \cup \{(A_v,B_v): v \in N_G(N_G[\{u\}]) \cap V(B_{\C_u})\}$.
Note that $\lvert N_G(N_G[\{u\}]) \cap V(B_{\C_u}) \rvert \leq \lvert V(A_{\C_u} \cap B_{\C_u}) \rvert \leq (d-1)\lvert \C_u \rvert \leq (d-1)d$.
So $\lvert \C_u' \rvert \leq \lvert \C_u \rvert + (d-1)d \leq d^2$.
Note that $N_G[N_G[\{u\}]] \cap V(B_{\C_u'})=\emptyset$.

Let $\C':=\bigcup_{z \in Z'} \C_{z}'$.
Then $N_G[N_G[Z']] \cap V(B_{\C'})=\emptyset$.
And $\lvert \C' \rvert \leq \lvert Z' \rvert d^2$.
Since $\lvert Y_1-Z \rvert \geq \lvert Y_1 \rvert - \lvert Z \rvert \geq 9\theta-\xi \geq 8\theta > \lvert Z \rvert$, there exists a subset $Y$ of $Y_1-Z$ with $\lvert Y \rvert = \lvert Z \rvert+1$.
Let $\C := \C' \cup \{(A_y,B_y): y \in Y\}$.
Clearly $(\C,\emptyset,\emptyset)$ satisfies (U1) and (U4).
Since $B_{\C'} \supseteq B_\C$, $(\C,\emptyset,\emptyset)$ satisfies (U2).
Since $Y \subseteq V(A_\C)-V(B_\C)$, $\lvert Z \rvert + 1 = \lvert Y \rvert \leq \lvert Y_1 \cap V(A_\C)-V(B_\C) \rvert$, so $(\C,\emptyset,\emptyset)$ satisfies (U3).
Note that $\lvert C \rvert \leq \lvert \C' \rvert + \lvert Y \rvert \leq \lvert Z'\rvert d^2 + \lvert Z \rvert +1$.
\end{proof}

For a useful triple $(\C,S,T)$, a vertex $v$ of $V(G)-Z$ is: 
	\begin{itemize}
		\item {\it $(\C,S,T)$-dangerous} if $v \in V(A_\C \cap B_\C)-S$ and there exists $v' \in N_G(v) \cap V(A_\C)-(V(B_\C) \cup S)$ such that either: 
			\begin{itemize}
				\item $v' \not \in Y_1$ and $\lvert ((Y_1 \cap V(A_\C)) \cup (S-V(A_\C \cap B_\C))) \cap N_G(v') \rvert \geq 2d-4$, or 
				\item $v' \in Y_1-T$,
			\end{itemize}
		\item {\it $(\C,S,T)$-heavy} if $v \in V(A_\C \cap B_\C)-S$ and 
		$$\lvert N_G(v) \cap ((Y_1 \cap V(A_\C)-V(B_\C)) \cup (S-V(A_\C \cap B_\C))) \rvert \geq d-1.$$
	\end{itemize}

\begin{claim}
\label{no H subdivision Kst subgraph Claim7}
 Let $(\C,S,T)$ be a useful triple and let $x \in V(A_\C \cap B_\C)$ be a $(\C,S,T)$-heavy vertex.
Then there exists a useful triple $(\C',S',T')$ with $\C'=\C \cup \{(A_x,B_x)\}$,  such that:
	\begin{itemize}
		\item $V(A_{\C'} \cap B_{\C'}) -S' \subseteq V(A_\C \cap B_\C)-S$,
		\item the set of $(\C',S',T')$-heavy vertices is strictly contained in the set of $(\C,S,T)$-heavy vertices, and 
		\item the set of $(\C',S',T')$-dangerous vertices is a subset of the set of $(\C,S,T)$-dangerous vertices.
	\end{itemize}
\end{claim}

\begin{proof}
Let $\C':=\C \cup \{(A_x,B_x)\}$.
Let $X:=N_G(x) \cap (Y_1 \cup S) \cap V(A_x \cap A_\C)-V(B_x \cup B_\C)$.
Since $x \in V(A_x)-V(B_x)$, $N_G(x) \subseteq V(A_x)$.
So $\lvert X \rvert \geq \lvert N_G(x) \cap (Y_1 \cup S) \cap V(A_\C)-V(B_\C) \rvert - \lvert V(A_x \cap B_x \cap A_\C)-V(B_\C) \rvert \geq d-1 - \lvert V(A_x \cap B_x \cap A_\C)-V(B_\C) \rvert$ since $x$ is $(\C,S,T)$-heavy.
That is, $\lvert V(A_x \cap B_x)-V(B_\C) \rvert = \lvert V(A_x \cap B_x \cap A_\C)-V(B_\C) \rvert \geq d-1 -\lvert X \rvert$.
Since $x$ is $(\C,S,T)$-heavy, $x \not \in S$.
Let $\iota$ be a bijection mentioned in (U3) witnessing that $(\C,S,T)$ is useful.
Let $X'$ be the intersection of $X$ and the domain of $\iota$.

For each $y \in X'$, since $x \in N_G(y) \cap V(A_\C \cap B_\C)-S$, $\iota(y) \in N_G(y) \cap V(A_\C \cap B_\C)$ by (U3).
Since $X' \subseteq X \subseteq V(A_x \cap A_\C)-V(B_x \cup B_\C)$, for each $y \in X'$, $N_G(y) \subseteq V(A_x)$, so $\iota(y) \in N_G(y) \cap V(A_\C \cap B_\C) \subseteq V(A_x) \cap V(A_\C \cap B_\C) \subseteq (V(A_{\C'})-V(B_{\C'})) \cup V(A_{\C} \cap B_{\C} \cap A_x \cap B_x)$.
Let 
\begin{align*}
Z_1 & := V(A_x \cap B_x)-V(A_\C)\\
Z_2 & := V(A_x \cap B_x \cap A_\C \cap B_\C)-\{\iota(y): y \in X'\} \\
Z_3 & :=V(A_x \cap B_x \cap A_\C \cap B_\C) \cap \{ \iota(y): y \in X' \}.
\end{align*}
So $\{Z_1, Z_2, Z_3\}$ is a partition of $V(A_x \cap B_x \cap B_\C)$, and hence 
\begin{align*}
\lvert Z_1 \cup Z_2 \cup Z_3 \rvert 
& = \lvert V(A_x \cap B_x \cap B_\C) \rvert\\ 
& = \lvert V(A_x \cap B_x) \rvert - \lvert V(A_x \cap B_x)-V(B_\C) \rvert \\
& \leq d-1- \lvert V(A_x \cap B_x)-V(B_\C) \rvert.
\end{align*}
Recall that $\lvert V(A_x \cap B_x)-V(B_\C) \rvert \geq d-1 -\lvert X \rvert$.
So $\lvert Z_1 \cup Z_2 \cup Z_3 \rvert \leq (d-1)-(d-1-\lvert X \rvert) = \lvert X \rvert$. 

Let $$S' := \big(S \cap V(B_x)-V(A_x \cap B_x \cap A_\C \cap B_\C)\big) \cup 
\big(V(A_{\C'} \cap B_{\C'})-V(A_\C \cap B_\C)\big) \cup 
V(A_x \cap B_x \cap A_\C \cap B_\C).$$
Note that $(V(A_{\C'} \cap B_{\C'})-V(A_\C \cap B_\C) \cup V(A_x \cap B_x \cap A_\C \cap B_\C) \subseteq V(A_x \cap B_x \cap B_\C) = Z_1 \cup Z_2 \cup Z_3$.
So 
\begin{align*}
\lvert S' \rvert 
& \leq (\lvert S \cap V(B_x) \rvert - \lvert S \cap V(A_x \cap B_x \cap A_\C \cap B_\C) \rvert) + \lvert Z_1 \cup Z_2 \cup Z_3 \rvert \\
& \leq \lvert S \cap V(B_x) \rvert - \lvert \{y \in X': \iota(y) \in Z_3\} \rvert + \lvert X \rvert \\
& = \lvert S \cap V(B_x) \rvert + \lvert X - \{y \in X': \iota(y) \in Z_3\} \rvert \\
& = \lvert \{y \in Y_1 \cap V(A_\C): \iota(y) \in S \cap V(B_x)\} \rvert + \lvert X-X' \rvert + \lvert X' - \{y \in X': \iota(y) \in Z_3\} \rvert.
\end{align*}
Recall that for every $y \in X'$, $\iota(y) \in N_G(y) \cap V(A_\C \cap B_\C) \cap V(A_x)$.
So if $y \in X'-\{y \in X': \iota(y) \in Z_3\}$, then 
$$\iota(y) \in N_G(y) \cap V(A_\C \cap B_\C) \cap V(A_x) - V(A_x \cap B_x \cap A_\C \cap B_\C) = N_G(y) \cap 
V(A_\C \cap B_\C \cap A_x) - V(B_x),$$ 
so $\iota(y) \not \in S \cap V(B_x)$. 
That is, $\{y \in Y_1 \cap V(A_\C): \iota(y) \in S \cap V(B_x)\}$ and $X' - \{y \in X': \iota(y) \in Z_3\}$ are disjoint.
Note that $X-X'$ is disjoint from the domain of $\iota$.
So $\{y \in Y_1 \cap V(A_\C): \iota(y) \in S \cap V(B_x)\}$, $X-X'$ and $X' - \{y \in X': \iota(y) \in Z_3\}$ are pairwise disjoint sets.
Therefore, $$\lvert S' \rvert \leq \lvert \{y \in Y_1 \cap V(A_\C): \iota(y) \in S \cap V(B_x)\} \cup (X - \{y \in X': \iota(y) \in Z_3\}) \rvert.$$

Since $X \subseteq Y_1 \cup S$ and $X \cap V(B_x)=\emptyset$, for every $x \in X - \{y \in X': \iota(y) \in Z_3\}$, if $x \not \in X \cap Y_1-\{y \in X': \iota(y) \in Z_3\}$, then $x \in X \cap S-Y_1$ and $\iota(y)=x$ for some $y \in Y_1 \cap V(A_\C)$ such that $\iota(y) \not \in S \cap V(B_x)$.
In addition, if $y$ is a vertex in $Y_1 \cap V(A_\C)$ such that $\iota(y) \in X \cap S-Y_1$, then $\iota(y) \not \in S \cap V(B_x)$.

Since $\lvert S' \rvert \leq \lvert \{y \in Y_1 \cap V(A_\C): \iota(y) \in S \cap V(B_x)\} \cup (X - \{y \in X': \iota(y) \in Z_3\}) \rvert$, there exists an injection $\iota'$ such that
	\begin{itemize}
		\item $\iota'(y)=\iota(y)$ if $y$ is in the domain of $\iota$ and $\iota(y) \in S \cap V(B_x)$,  
		\item for each $v \in S'-(S \cap V(B_x))$, there exists exactly one element $y \in (X \cap Y_1-\{y \in X': \iota(y) \in Z_3\}) \cup \{y \in Y_1 \cap V(A_\C): \iota(y) \in X \cap S-Y_1\}$ such that $\iota'(y)=v$, and
		\item if $\iota(y_1)=\iota'(y_2)$ for some $y_1,y_2$, then $y_1=y_2$.
	\end{itemize}
Recall that $\iota(y) \not \in S \cap V(B_x)$ for every $y \in (X \cap Y_1-\{y \in X': \iota(y) \in Z_3\}) \cup \{y \in Y_1: \iota(y) \in X \cap S-Y_1\}$.
Then $\iota'$ is a bijection from a subset of $Y_1 \cap V(A_{\C'})$ to $S'$.
We further modify $\iota'$ and $S'$ by applying the following operations  for some vertex $y \in V(A_{\C'})-V(B_{\C'})$ in the domain of $\iota'$ with $\iota'(y) \not \in N_G(y) \cap V(A_{\C'} \cap B_{\C'})$ and $N_G(y) \cap V(A_{\C'} \cap B_{\C'})-S' \neq \emptyset$, and then repeating until no such vertex $y$ exists: 
	\begin{itemize}
		\item add a vertex $v \in N_G(y) \cap V(A_{\C'} \cap B_{\C'})-S'$ into $S'$,
		\item delete $\iota'(y)$ from $S'$, and 
		\item redefine $\iota'(y)$ to be $v$.''
	\end{itemize}
Now, further modify $\iota'$ and $S'$ by applying the following operations for some vertex $z \in S' - N_G[V(A_{\C'} \cap B_{\C'})]$, and repeating until no such vertex $z$ exists:
	\begin{itemize}
		\item remove $z$ from $S'$, and
		\item if $y$ is the element in the domain of $\iota'$ with $\iota'(y)=z$, then remove $y$ from the domain of $\iota'$.
	\end{itemize}	
	
Notice that for each vertex $z$ removed from $S'$ in the above procedure, $z \in S-V(A_\C \cap B_\C)$ and $N_G(z) \cap V(A_\C \cap B_\C) \subseteq V(A_\C \cap B_\C)-V(B_x)$.
Note that $\iota'$ remains a bijection from a subset of $Y_1 \cap V(A_{\C'})$ to $S'$.

Observe that for every $y$ in the domain of $\iota'$ with $y \in V(A_\C')-V(B_\C')$ and $N_G(y) \cap V(A_{\C'} \cap B_{\C'})-S' \neq \emptyset$, $\iota'(y) \in N_G(y) \cap V(A_{\C'} \cap B_{\C'})$ due to the above modification.
In addition, if $y$ is in the domain of $\iota'$ and $y \in V(A_{\C'} \cap B_{\C'})$, then $y \in V(A_\C \cap B_\C)$ and $y$ is in the domain of $\iota$ such that $\iota(y)=\iota'(y)$, so $\iota'(y)=\iota(y)=y$.

Let $T'$ be the set obtained from $T$ by deleting the domain of $\iota'$.
So $T'$ is disjoint from the domain of $\iota'$.
Since $T$ is disjoint from $Z'$, $T'$ is disjoint from $Z'$.
So $(\C',S',T')$ satisfies (U4).
In addition, $\lvert S' \rvert - \lvert S \rvert$ is at most the number of vertices in $X$ and in the domain of $\iota'$ but not in the domain of $\iota$.
So $\lvert S' \rvert - \lvert S \rvert \leq \lvert T \rvert - \lvert T' \rvert$.
Hence $\lvert S' \rvert + \lvert T' \rvert \leq \lvert S \rvert + \lvert T \rvert$.
Since $\{y \in Y_1 \cap S \cap V(A_\C \cap B_\C): \iota(y)=y\}-\{y \in Y_1 \cap S' \cap V(A_{\C'} \cap B_{\C'}): \iota'(y)=y\} \subseteq (V(A_{\C'})-V(B_{\C'})) - (V(A_\C)-V(B_\C))$, $(\C',S',T')$ satisfies (U3) and is useful.

It is easy to see that $V(A_{\C'} \cap B_{\C'})-S' \subseteq V(A_\C \cap B_\C)-S$.
Note that each vertex $v \in V(A_{\C'} \cap B_{\C'})-S'$ belongs to $V(A_\C \cap B_\C)-V(A_x)$, 
so $N_G(v) \cap V(A_{\C'})-V(B_{\C'})=N_G(v) \cap V(A_\C)-V(B_\C)$.
Furthermore, $S'-V(A_{\C'} \cap B_{\C'}) \subseteq S-V(A_\C \cap B_\C)$.
Hence every $(\C',S',T')$-heavy vertex is $(\C,S,T)$-heavy.
Since $x$ is $(\C,S,T)$-heavy but not $(\C',S',T')$-heavy, the set of $(\C',S',T')$-heavy vertices is strictly contained in the set of $(\C,S,T)$-heavy vertices.

Let $v$ be a $(\C',S',T')$-dangerous vertex and let $v'$ be a vertex in $N_G(v) \cap V(A_{\C'})-(V(B_{\C'}) \cap S')$ witnessing the definition of being dangerous.
Since $v \not \in S'$, $v \in V(A_\C \cap B_\C)-V(A_x)$, so $v' \in N_G[V(A_\C \cap B_\C)] \cap V(A_\C) \cap V(B_x)$.
So $v' \in V(B_x)-(V(B_\C) \cup S')$.
Since $v \in V(A_{\C'} \cap B_{\C'})$, $v' \in N_G(v) \cap N_G[V(A_{\C'} \cap B_{\C'})]$.
Since $v' \not \in S'$ and $v' \not \in V(A_\C \cap B_\C)-V(A_{\C'} \cap B_{\C'})$ and $v' \in V(B_x)$ and $N_G(v') \cap V(A_\C \cap B_\C)-V(A_x) \neq \emptyset$, we know $v' \not \in S$ by the procedure of modifying $S$.
So $v' \in (N_G(v) \cap V(A_\C)-(V(B_\C) \cup S)) \cap V(B_x)$.
Note that $T-T' \subseteq V(A_x)-V(B_x)$.
So if $v' \in Y_1-T'$, then $v' \in Y_1-T$ and $v$ is $(\C,S,T)$-dangerous.
Furthermore, $Y_1 \cap V(A_\C) \cap N_G(v') = Y_1 \cap V(A_{\C'}) \cap N_G(v')$ and $S' - V(A_{\C'} \cap B_{\C'}) \subseteq S -V(A_\C \cap B_\C)$, so $v$ is $(\C,S,T)$-dangerous.
Therefore, the set of $(\C',S',T')$-dangerous vertices is a subset of the set of $(\C,S,T)$-dangerous vertices.
This proves the claim.
\end{proof}

\begin{claim}
\label{no H subdivision Kst subgraph Claim8}
Let $(\C,S,T)$ be a useful triple.
Then there exists a set $S'$ with $S \cup (Y_1 \cap V(A_\C \cap B_\C)) \subseteq S' \subseteq N_G[V(A_\C \cap B_\C)] \cap V(A_\C)$ such that $(\C,S',T)$ is a useful triple and:
	\begin{itemize}
		\item If $\iota'$ is the bijection witnessing that $(\C,S',T)$ satisfies (U3), then for every $y \in Y_1 \cap V(A_\C \cap B_\C)$, the unique element of the domain of $\iota'$ mapped to $y$ by $\iota'$ is $y$. 
		\item The set of $(\C,S',T)$-dangerous vertices is contained in the set of $(\C,S,T)$-dangerous vertices.
		\item The set of $(\C,S',T)$-heavy vertices is contained in the set of $(\C,S,T)$-heavy vertices.
	\end{itemize}
\end{claim}

\begin{proof}
Let $\iota$ be a function mentioned in (U3) witnessing that $(\C,S,T)$ is a useful triple.
We may assume that $Y_1 \cap V(A_\C \cap B_\C) \subseteq S$, since if some vertex $y \in Y_1 \cap V(A_\C \cap B_\C)$ does not belong to $S$, then $y$ is not in the domain of $\iota$, and we can define $\iota(y)=y$ without violating (U3) and (U4) such that the set of dangerous vertices and the set of heavy vertices remain the same.

Since $\iota$ is a bijection, we write the element mapped to $y$ by $\iota$ as $\iota^{(-1)}(y)$.
Modify $\iota$ and $S$ by applying the following operations to some vertex $y \in Y_1 \cap S \cap V(A_\C \cap B_\C)$ with $\iota^{(-1)}(y) \neq y$,  and repeat until no such $y$ exists:
	\begin{itemize}
		\item remove $\iota^{(-1)}(y)$ from the domain of $\iota$, 
		\item define $\iota(y):=y$,
	\end{itemize}

Then define $S'$ and $\iota'$ to be the modified $S$ and $\iota$, respectively.
Clearly, $(\C,S',T)$ satisfies (U3), $S \subseteq S' \subseteq N_G[V(A_\C \cap B_\C)] \cap V(A_\C)$, and $\iota'(y)=y$ for every $y \in Y_1 \cap S' \cap V(A_\C \cap B_\C)$.
Since we assume that $Y_1 \cap V(A_\C \cap B_\C) \subseteq S$, we have $S \cup (Y_1 \cap V(A_\C \cap B_\C)) \subseteq S' \subseteq N_G[V(A_\C \cap B_\C)] \cap V(A_\C)$ and $\iota'(y)=y$ for every $y \in Y_1 \cap V(A_\C \cap B_\C)$.
Since $T \subseteq V(A_\C)-V(B_\C)$, $(\C,S',T)$ satisfies (U4).
Since $S'-S \subseteq V(A_\C \cap B_\C)$, the set of $(\C,S',T)$-dangerous vertices is contained in the set of $(\C,S,T)$-dangerous vertices, and the set of $(\C,S',T)$-heavy vertices is contained in the set of $(\C,S,T)$-heavy vertices.
\end{proof}

\begin{claim}
\label{no H subdivision Kst subgraph Claim9}
Let $(\C,S,T)$ be a useful triple, and let $x$ be a $(\C,S,T)$-dangerous vertex.
If there exists no $(\C,S,T)$-heavy vertex, then there exists a useful triple $(\C',S',T')$ with $\C'=\C \cup \{(A_x,B_x)\}$ such that the set of $(\C',S',T')$-dangerous vertices is strictly contained in the set of $(\C,S,T)$-dangerous vertices.
\end{claim}

\begin{proof}
By \cref{no H subdivision Kst subgraph Claim8}, we may assume that $Y_1 \cap V(A_\C \cap B_\C) \subseteq S$ and the function $\iota$ mentioned in (U3) witnessing that $(\C,S,T)$ is useful satisfies $\iota(y)=y$ for every $y \in Y_1 \cap V(A_\C \cap B_\C)$. Let $\C':=\C \cup \{(A_x,B_x)\}$.

We first assume that $\lvert Y_1 \cap V(A_x \cap B_\C) \rvert \geq d-2$.
So $\lvert Y_1 \cap V(A_x \cap B_\C) \rvert \geq \lvert V(A_x \cap B_x) \cap V(B_\C) \rvert$ by \cref{no H subdivision Kst subgraph Claim5}.
Hence there exists a function $\iota'$ whose domain is a subset of $Y_1 \cap V(A_\C \cup A_x)$ such that:
	\begin{itemize}
		\item $\iota'(y)=\iota(y)$ for every $y \in Y_1 \cap V(A_\C) -V(A_x \cap B_\C)$ belonging to the domain of $\iota$ with $\iota(y) \in V(A_\C) \cap N_G[V(A_{\C'} \cap B_{\C'})]-V(A_x)$, and 
		\item for each vertex $v$ in $V(A_x \cap B_x) \cap V(B_\C)$, there exists exactly one element $y \in Y_1 \cap V(A_x \cap B_\C)$ such that $\iota'(y)=v$ and if $v \in Y_1$, then $y=v$.
	\end{itemize}
Let $S':=(S \cap N_G[V(A_{\C'} \cap B_{\C'})] - V(A_x)) \cup V(A_x \cap B_x \cap B_\C)$.
So $\iota'$ is a bijection from a subset of $Y_1 \cap V(A_\C \cup A_x)$ to $S'$.
Note that every vertex in $S'-V(A_{\C'} \cap B_{\C'})$ is contained in $S -(V(A_x) \cup V(A_{\C'} \cap B_{\C'}))$, so it is adjacent to some vertex in $V(A_{\C'} \cap B_{\C'})$.

Let $T':=T$.
Since $\iota$ satisfies (U3) and $\iota(y)=y$ for every $y \in Y_1 \cap S \cap V(A_\C \cap B_\C)$, we know that $\iota'$ satisfies (U3).
Since $T'=T \subseteq V(A_\C)-V(B_\C)$, $(\C',S',T')$ satisfies (U4).
So $(\C',S',T')$ is a useful triple.

Let $v$ be a $(\C',S',T')$-dangerous vertex.
So $v \in V(A_{\C'} \cap B_{\C'})-S' \subseteq V(A_\C \cap B_\C)-V(A_x)$.
Let $v'$ be a vertex witnessing that $v$ is $(\C',S',T')$-dangerous.
So $v' \in V(B_x)-(V(B_\C) \cup S')$ and $N_G(v') \subseteq V(A_\C)$.
Since $v \in V(A_{\C'} \cap B_{\C'})-V(A_x)$, $v' \in N_G[V(A_{\C'} \cap B_{\C'})]$, so $v' \not \in S$.
Since $S'-V(A_{\C'} \cap B_{\C'}) \subseteq S-V(A_\C \cap B_\C)$, if $v' \not \in Y_1$ and $\lvert ((Y_1 \cap V(A_{\C'})) \cup (S'-V(A_{\C'} \cap B_{\C'}))) \cap N_G(v') \rvert \geq 2d-4$, then $v' \not \in Y_1$ and $\lvert ((Y_1 \cap V(A_\C)) \cup (S-V(A_\C \cap B_\C))) \cap N_G(v') \rvert \geq 2d-4$, so $v$ is $(\C,S,T)$-dangerous.
Since $T'=T$, if $v' \in Y_1-T'$, then $v' \in Y_1-T$ and $v$ is $(\C,S,T)$-dangerous.
So every $(\C',S',T')$-dangerous vertex is $(\C,S,T)$-dangerous. 
Since $x$ is $(\C,S,T)$-dangerous but not $(\C',S',T')$-dangerous, the set of $(\C',S',T')$-dangerous vertices is strictly contained in the set of $(\C,S,T)$-dangerous vertices.
So the claim holds.

Hence we may assume that $\lvert Y_1 \cap V(A_x \cap B_\C) \rvert \leq d-3$.

Modify $S$ and define $\iota'$ to be the function obtained from $\iota$ by applying the following operations to a vertex $y$ in the domain of $\iota$ with $\iota(y) \not \in N_G[V(A_{\C'} \cap B_{\C'})] \cap V(A_{\C'})$, and  repeating until no such $y$ exists:
	\begin{itemize}
		\item if $y \in V(A_\C \cap B_\C)-V(B_x)$ or $V(A_x \cap B_x) \cap V(B_\C)-S = \emptyset$, then remove $y$ from the domain of $\iota$ and remove $\iota(y)$ from $S$, 
		\item if $y \not\in V(A_\C \cap B_\C)-V(B_x)$ and $N_G(y) \cap V(A_{\C'} \cap B_{\C'})-S=\emptyset$ and $V(A_x \cap B_x) \cap V(B_\C)-S \neq \emptyset$, then redefine $\iota(y)$ to be an element in $V(A_x \cap B_x) \cap V(B_\C)-S$ and add this element into $S$,
		\item otherwise remove $\iota(y)$ from $S$, redefine $\iota(y)$ to be an element in $N_G(y) \cap V(A_{\C'} \cap B_{\C'})-S$ and add this element into $S$.
	\end{itemize}
Let $S'$ be the modified $S$, and let 
$$T':=T \cup (Y_1 \cap V(B_\C)-V(A_\C \cup B_x)) \cup (Y_1 \cap V(A_\C \cap B_\C)-(V(B_x) \cup N_G[V(A_{\C'} \cap B_{\C'})])).$$

Clearly, $T'$ is disjoint from the domain of $\iota'$.
By (U2), $Z'  \cap V(B_\C)=\emptyset$.
So $T'$ is disjoint from $Z'$ as $T$ is disjoint from $Z'$.
So $(\C',S',T')$ satisfies (U4).

Since $Y_1 \cap V(A_\C \cap B_\C) \subseteq S$, we know 
$Y_1 \cap V(A_\C \cap B_\C)-(V(B_x) \cup N_G[V(A_{\C'} \cap B_{\C'})]) \subseteq S$.
For every $y \in Y_1 \cap V(A_\C \cap B_\C)-(V(B_x) \cup N_G[V(A_{\C'} \cap B_{\C'})])$, since $\iota(y)=y \not \in N_G[V(A_{\C'} \cap B_{\C'})] \cap V(A_{\C'})$ and $y \in V(A_\C \cap B_\C)-V(B_x)$, $y \in S-S'$.

So $\lvert S \rvert \geq \lvert S' \rvert + \lvert Y_1 \cap V(A_\C \cap B_\C)-(V(B_x) \cup N_G[V(A_{\C'} \cap B_{\C'})]) \rvert$.
Hence $\lvert S' \rvert + \lvert T' \rvert \leq \lvert S \rvert + \lvert T \rvert+\lvert Y_1 \cap V(B_\C)-V(A_\C \cup B_x) \rvert$.
Since $Y_1 \cap V(B_\C)-V(A_\C \cup B_x) \subseteq Y_1 \cap (V(A_{\C'})-V(B_{\C'}))-V(A_\C)$ and $(\C,S,T)$ satisfies (U3), we know 
\begin{align*}
& \;\;\;\;\;  \lvert S' \rvert + \lvert T' \rvert + \lvert Z \rvert + 1  \\
& \leq   \lvert S \rvert + \lvert T \rvert +\lvert Y_1 \cap V(B_\C)-V(A_\C \cup B_x) \rvert+ \lvert Z \rvert + 1 \\
& \leq  \lvert Y_1 \cap V(A_\C) - V(B_\C) \rvert + \lvert \{y \in Y_1 \cap S \cap V(A_\C \cap B_\C): \iota(y) =y\} \rvert \\ 
 & \quad \quad + \lvert Y_1 \cap (V(A_{\C'})-V(B_{\C'}))-V(A_\C) \rvert  \\
& \leq   \lvert Y_1 \cap V(A_{\C'})-V(B_{\C'}) \rvert - \lvert Y_1 \cap V(A_\C \cap B_\C) - V(B_x) \rvert \\
 &\quad \quad + \lvert \{y \in Y_1 \cap S \cap V(A_\C \cap B_\C): \iota(y) =y\} \rvert \\
& \leq \lvert Y_1 \cap V(A_{\C'})-V(B_{\C'}) \rvert + \lvert \{y \in Y_1 \cap S \cap V(A_\C \cap B_\C) \cap V(B_x): \iota(y)=y\} \rvert \\
& \leq  \lvert Y_1 \cap V(A_{\C'})-V(B_{\C'}) \rvert + \lvert \{y \in Y_1 \cap S' \cap V(A_{\C'} \cap B_{\C'}): \iota'(y)=y\} \rvert.
\end{align*}
Hence $(\C',S',T')$ satisfies (U3).
Therefore $(\C',S',T')$ is useful.

Suppose that the set of $(\C',S',T')$-dangerous vertices is not strictly contained in the set of $(\C,S,T)$-dangerous vertices.
Since $x$ is $(\C,S,T)$-dangerous but not $(\C',S',T')$-dangerous, there exists a vertex $v$ that is $(\C',S',T')$-dangerous but not $(\C,S,T)$-dangerous.
So there exists a vertex $v' \in N_G(v) \cap V(A_{\C'})-(V(B_{\C'}) \cup S')$ such that either $v' \in Y_1-T'$, or $v' \not \in Y_1$ and $\lvert ((Y_1 \cap V(A_{\C'})) \cup (S'-V(A_{\C'} \cap B_{\C'})) \cap N_G(v') \rvert \geq 2d-4$.
Since $v' \in N_G[V(A_{\C'} \cap B_{\C'})] \cap V(A_{\C'})$, if $v'$ belongs to $S$ at beginning, then $v'$ is not removed from $S$ during the process of modifying $S$, so $v' \in S'$, a contradiction.
So $v' \not \in S$.

Suppose that $v' \in V(A_\C)-V(B_\C)$.
So $v \in V(A_\C \cap B_\C) \cap V(A_{\C'} \cap B_{\C'}) \subseteq V(A_\C \cap B_\C) \cap V(B_x)$.
Hence if $v$ belongs to $S$ at beginning, then $v$ is not removed from $S$ during the process of modifying $S$, so $v \in S'$.
Since $v$ is $(\C',S',T')$-dangerous, $v \not \in S'$, so $v \not \in S$.
Since $v$ is not $(\C,S,T)$-dangerous, $v' \not \in Y_1-T$, and either $v' \in Y_1$ or $\lvert ((Y_1 \cap V(A_{\C})) \cup (S-V(A_\C \cap B_\C))) \cap N_G(v') \rvert < 2d-4$.
Since $T' \cap V(A_\C)-V(B_\C) = T \cap V(A_\C)-V(B_\C)$, $v' \not \in Y_1-T'$.
Since $v$ is $(\C',S',T')$-dangerous, $v' \not \in Y_1$ and $\lvert ((Y_1 \cap V(A_{\C'})) \cup (S'-V(A_{\C'} \cap B_{\C'}))) \cap N_G(v') \rvert \geq 2d-4$.
Since $N_G(v') \subseteq V(A_\C)$ and $S'-V(A_{\C'} \cap B_{\C'}) \subseteq S-V(A_\C \cap B_\C)$, 
\begin{align*}
2d-4 & \leq \lvert ((Y_1 \cap V(A_{\C'})) \cup (S'-V(A_{\C'} \cap B_{\C'}))) \cap N_G(v') \rvert\\ 
& = \lvert ((Y_1 \cap V(A_\C)) \cup (S'-V(A_{\C'} \cap B_{\C'}))) \cap N_G(v') \rvert\\ 
& \leq \lvert ((Y_1 \cap V(A_\C)) \cup (S-V(A_{\C} \cap B_{\C}))) \cap N_G(v') \rvert\\
& < 2d-4,\end{align*}
a contradiction.

Therefore, $v' \in V(B_\C)$.
So $v' \in (V(A_x)-V(B_x)) \cap V(B_\C)$ and hence $N_G(v') \subseteq V(A_x)$.
Since $Y_1 \cap V(A_\C \cap B_\C) \subseteq S$, if $v' \in Y_1 \cap V(A_\C \cap B_\C)$, then $v' \in S$, a contradiction.
So $v' \not \in Y_1 \cap V(A_\C \cap B_\C)$.
Since $Y_1 \cap V(B_\C)-V(A_\C \cup B_x) \subseteq T'$, if $v' \in Y_1-T'$, then $v' \in Y_1 \cap V(A_\C \cap B_\C)$, a contradiction.
So $v' \not \in Y_1-T'$.
Since $v$ is $(\C',S',T')$-dangerous, $\lvert ((Y_1 \cap V(A_{\C'})) \cup (S'-V(A_{\C'} \cap B_{\C'}))) \cap N_G(v') \rvert \geq 2d-4$.
Since $\lvert Y_1 \cap V(A_x \cap B_\C) \rvert \leq d-3$ and $(S'-V(A_{\C'} \cap B_{\C'}))-(V(A_\C) \cup V(B_x))=\emptyset$ and $N_G(v') \subseteq V(A_x)$, 
\begin{align*}
 &\; \lvert N_G(v') \cap ((Y_1 \cap V(A_\C)-V(B_\C)) \cup (S-V(A_{\C} \cap B_{\C})) \rvert \\
\geq &\; \lvert N_G(v') \cap ((Y_1 \cap V(A_x)-V(B_\C)) \cup (S'-V(A_{\C'} \cap B_{\C'})) \rvert \\
\geq &\; \lvert N_G(v') \cap ((Y_1 \cap V(A_x)) \cup (S'-(V(A_{\C'} \cap B_{\C'})))) \rvert - (d-3)\\
= &\; \lvert N_G(v') \cap (Y_1 \cup (S'-(V(A_{\C'} \cap B_{\C'})))) \rvert - (d-3)\\
 \geq &\;  (2d-4)-(d-3) \\
 = &\; d-1.
\end{align*}
In particular, $v' \in V(A_\C \cap B_\C)-V(B_x)$.
Since there exists no $(\C,S,T)$-heavy vertex, $v'$ is not a $(\C,S,T)$-heavy vertex.
So $v' \in S$, a contradiction.
This proves the claim.
\end{proof}

\begin{claim}
\label{no H subdivision Kst subgraph Claim10}
If $(\C,S,T)$ is a useful triple such that there exists a $(\C,S,T)$-dangerous vertex, then there exists a useful triple $(\C',S',T')$ with $\C \subseteq \C'$ and $\lvert \C' \rvert \leq \lvert \C \rvert + \lvert V(A_\C \cap B_\C) \rvert+1$ such that the set of $(\C',S',T')$-dangerous vertices is strictly contained in the set of $(\C,S,T)$-dangerous vertices.
\end{claim}

\begin{proof}
Note that there are at most $\lvert V(A_\C \cap B_\C) \rvert$ $(\C,S,T)$-heavy vertices.
By repeatedly applying \cref{no H subdivision Kst subgraph Claim7} at most $\lvert V(A_\C \cap B_\C) \rvert$ times, there exists a useful triple $(\C_1,S_1,T_1)$ with $\C \subseteq \C_1$ and $\lvert \C_1 \rvert \leq \lvert \C \rvert + \lvert V(A_\C \cap B_\C) \rvert$ such that there exists no $(\C_1,S_1,T_1)$-heavy vertices, and the set of $(\C_1,S_1,T_1)$-dangerous vertices is contained in the set of $(\C,S,T)$-dangerous vertices.
By \cref{no H subdivision Kst subgraph Claim9}  applied to $\C_1$, there exists a useful triple $(\C',S',T')$ with $\C_1 \subseteq \C'$ and $\lvert \C' \rvert=\lvert \C_1 \rvert+1 \leq \lvert \C \rvert + \lvert V(A_\C \cap B_\C) \rvert+1$ such that the set of $(\C',S',T')$-dangerous vertices is strictly contained in the set of $(\C_1,S_1,T_1)$-dangerous vertices and hence is strictly contained in the set of $(\C,S,T)$-dangerous vertices.
This proves the claim.
\end{proof}

\begin{claim}
\label{no H subdivision Kst subgraph Claim11}
There exists a useful triple $(\C^*,S^*,T^*)$ with $\lvert \C^* \rvert \leq a_{(d-1)a_0}$ such that there exists no $(\C^*,S^*,T^*)$-dangerous vertex.
\end{claim}

\begin{proof}
By \cref{no H subdivision Kst subgraph Claim6}, there exists a useful triple $(\C_0,\emptyset,\emptyset)$ with $\lvert \C_0 \rvert \leq \lvert Z' \rvert d^2 + \lvert Z \rvert+1$.
Let $S_0=\emptyset$ and $T_0=\emptyset$.
So $(\C_0,S_0,T_0)$ is a useful triple with $\lvert \C_0 \rvert \leq f(\xi)d^2 + \xi+1=a_0$.
For $i \geq 1$, if there exists a $(\C_{i-1},S_{i-1},T_{i-1})$-dangerous vertex, then by 
\cref{no H subdivision Kst subgraph Claim10}, there exists a useful triple $(\C_i,S_i,T_i)$ such that $\lvert \C_i \rvert \leq \lvert \C_{i-1} \rvert + \lvert V(A_{\C_{i-1}} \cap B_{\C_{i-1}}) \rvert+1$ and the set of $(\C_i,S_i,T_i)$-dangerous vertices is strictly contained in the set of $(\C_{i-1},S_{i-1},T_{i-1})$-dangerous vertices. 
So $\lvert \C_i \rvert \leq a_{i-1}+(d-1)a_{i-1}+1 \leq a_i$ for each $i \geq 1$ by induction on $i$.
Since there are at most $\lvert V(A_{\C_0} \cap B_{\C_0}) \rvert \leq \lvert \C_0 \rvert(d-1) \leq (d-1)a_0$ $(\C_0,S_0,T_0)$-dangerous vertices.
Hence there exists $i^*$ with $0 \leq i^* \leq (d-1)a_0$ such that $(\C_{i^*},S_{i^*},T_{i^*})$ is a useful triple with no $(\C_{i^*},S_{i^*},T_{i^*})$-dangerous vertex.
Note that $\lvert \C_{i^*} \rvert \leq a_{i^*} \leq a_{(d-1)a_0}$. 
\end{proof}

Let $\iota^*$ be the function mentioned in (U3) witnessing that $(\C^*,S^*,T^*)$ is useful.
By \cref{no H subdivision Kst subgraph Claim8}, we may assume that $Y_1 \cap V(A_{\C^*} \cap B_{\C^*}) \subseteq S^*$ such that $\iota^*(y)=y$ for every $y \in Y_1 \cap V(A_{\C^*} \cap B_{\C*})$.

Define the following:
	\begin{align*}
G_B & :=G[(\bigcap_{(A,B) \in \C^*} V(B)) \cup (Z \cup S^* \cup T^*)]\\
Y_B & :=(Y_1 \cap V(B_{\C^*})) \cup (Z \cup S^* \cup T^*).
	\end{align*}

\begin{claim} \label{claim:Y_B}
For every vertex $v \in V(G_B)-Y_B$, $N_G(v) \cap Y_1 \subseteq N_{G_B}(v) \cap Y_B$.
\end{claim}

\begin{proof}
Suppose to the contrary that there exist $v \in V(G_B)-Y_B$ and $y \in N_G(v) \cap Y_1-(N_{G_B}(v) \cap Y_B)$.
Since $y \in Y_1-Y_B$, $y \in V(A_{\C^*})-V(B_{\C^*})$.
So $v \in V(A_{\C^*} \cap B_{\C^*})-S^*$.
Since there exists no $(\C^*,S^*,T^*)$-dangerous vertex, $v$ is not a $(\C^*,S^*,T^*)$-dangerous vertex.
Since $y \in N_G(v) \cap V(A_{\C^*})-(V(B_{\C^*}) \cup S^*)$, $y \not \in Y_1-T^*$.
Since $y \in Y_1$, $y \in T^*$.
So $y \in Y_B$, a contradiction.
\end{proof}

Define the following:
	\begin{itemize}
		\item For every $y \in Y_B$, let $L_B(y)$  be a 1-element subset of $L(y)$.
		\item For every $v \in V(G_B)-Y_B$ with $\lvert N_{G_B}(v) \cap Y_B \rvert \in [ s'-1]$, let $L_B(v)$ be a subset of $L(v)$ with size $s'+2-\lvert N_{G_B}(v) \cap Y_B \rvert$ such that $L_B(v) \cap L_B(u)=\emptyset$ for every $u \in N_G(v) \cap Y_B$.
			(Note that such a subset of $L(v)$ exists by \cref{claim:Y_B}.)
		\item For every other vertex $v$ of $G_B$, let $L_B(v):=L(v)$.
	\end{itemize}
Hence $L_B$ is a restricted $(s',2,Y_B)$-list-assignment by \cref{claim:Y_B}.
Since $(\C^*,S^*,T^*)$ is useful and $\iota^*(y)=y$ for every $y \in Y_1 \cap V(A_{\C^*} \cap B_{\C*})$, 
\begin{align*}
\lvert Y_B \rvert 
&  \leq \lvert Y_1 \cap V(B_{\C^*}) \rvert + \lvert S^* \rvert - \lvert Y_1 \cap V(A_{\C^*} \cap B_{\C^*}) \rvert + \lvert Z \rvert + \lvert T^* \rvert \\
& \leq (\lvert Y_1 \rvert - \lvert Y_1 \cap V(A_{\C^*})-V(B_{\C^*}) \rvert) + \lvert S^* \rvert-\lvert Y_1 \cap V(A_{\C^*} \cap B_{\C^*}) \rvert + \lvert Z \rvert + \lvert T^* \rvert \\
& \leq \lvert Y_1 \rvert-1
\end{align*} 
by (U3).

Since $\lvert Y_B \rvert < \lvert Y_1 \rvert$, we know $\lvert V(G_B) \rvert < \lvert V(G) \rvert$.
By the minimality of $\lvert V(G) \rvert$, there exists an $(\eta,g)$-bounded $L_B$-coloring $c_B$ of $G_B$.
Define the following:
	\begin{itemize}
		\item Let $G_A:=G[V(A_{\C^*}) \cup Z]$.
		\item Let $Y_A:=(Y_1 \cap V(A_{\C^*})) \cup Z \cup S^* \cup T^* \cup V(A_{\C^*} \cap B_{\C^*})$.
		\item For every $y \in Y_A$, let $L_A(y)$ be a 1-element subset of $L(y)$ such that if $y \in V(G_B)$, then $L_A(y)=\{c_B(y)\}$.
		\item For every $v \in V(G_A)-Y_A$ with $1 \leq \lvert N_{G_A}(v) \cap Y_A \rvert \leq s'-1$, 
		let $L_A(v)$ be a subset of $L(v)$ with size $s'+2-\lvert N_{G_A}(v) \cap Y_A \rvert$ such that $L_A(v) \cap L_A(u)=\emptyset$ for every $u \in Y_A \cap N_{G_A}(v)$.
		\item For every other vertex $v$ of $G_A$, let $L_A(v) := L(v)$.
	\end{itemize}
Then $L_A$ is a restricted $(s',2,Y_A)$-list-assignment of $G_A$.
Since $\lvert \C^* \rvert \leq a_{(d-1)a_0}$, $\lvert V(A_{\C^*} \cap B_{\C^*}) \rvert \leq (d-1)a_{(d-1)a_0}<\theta-\xi$.
So $(A_{\C^*},B_{\C^*}) \in \T-Z$ and hence $\lvert Y_1 \cap V(A_{\C^*}) \rvert \leq 3\theta$.
By (U3), $\lvert S^* \rvert + \lvert T^* \rvert + \lvert Z \rvert \leq \lvert Y_1 \cap V(A_{\C^*})-V(B_{\C^*}) \rvert + \lvert Y_1 \cap V(A_{\C^*} \cap B_{\C^*}) \rvert \leq \lvert Y_1 \cap V(A_{\C^*}) \rvert \leq 3\theta$.
Hence $\lvert Y_A \rvert \leq \lvert Y_1 \cap V(A_{\C^*}) \rvert + \lvert S^* \rvert + \lvert T^* \rvert + \lvert Z \rvert + \lvert V(A_{\C^*} \cap B_{\C^*}) \rvert \leq 3\theta + 3\theta+\theta \leq 7\theta$.

In particular, $\lvert V(G_A) \rvert < \lvert V(G) \rvert$.
By minimality, there exists an $(\eta,g)$-bounded $L_A$-coloring $c_A$ of $G_A$.

\begin{claim}
\label{no H subdivision Kst subgraph Claim12} 
For every $v \in V(G_A)-(V(G_B) \cup Y_1)$ with $N_G(v) \cap V(A_{\C^*} \cap B_{\C^*})-S^* \neq \emptyset$, $c_B(u) \not \in L_A(v)$ for every $u \in N_G(v) \cap V(G_B)$.
\end{claim}

\begin{proof}
Since $v \in V(G_A)-(V(G_B) \cup Y_1)$, $v \in V(G)-(Z \cup V(B_{\C^*}))$.
So there exists $(A,B) \in \C^*$ such that $v \in V(A)-V(B)$.
Hence $N_G(v) \subseteq V(A)$ and $\lvert N_G(v) \cap V(A_{\C^*} \cap B_{\C^*}) \rvert \leq \lvert N_G(v) \cap V(A \cap B) \rvert \leq d-1$.

Since $v \in N_G(A_{\C^*} \cap B_{\C^*})$ and $\C^*$ satisfies (U2), $v \not \in Z'$.
So $\lvert N_G(v) \cap Z \rvert \leq s-1$.

Since $N_G(v) \cap V(A_{\C^*} \cap B_{\C^*})-S^* \neq \emptyset$, there exists $w \in N_G(v) \cap V(A_{\C^*} \cap B_{\C^*})-S^*$.
Since there exists no $(\C^*,S^*,T^*)$-dangerous vertex, $w$ is not a $(\C^*,S^*,T^*)$-dangerous vertex.
Since $S^* \subseteq V(G_B)$, $v \not \in S^*$.
Since $v \not \in Y_1$, $v \not \in Y_1-T^*$.
So $\lvert N_G(v) \cap ((Y_1 \cap V(A_{\C^*})) \cup (S^*-V(A_{\C^*} \cap B_{\C^*}))) \rvert \leq 2d-5$.

Since $T^* \subseteq Y_1 \cap V(A_{\C^*})$, 
\begin{align*}
& \,\,\,\,\,\, \,\, \lvert N_G(v) \cap Y_A \rvert \\
& \leq \lvert N_G(v) \cap V(A_{\C^*} \cap B_{\C^*}) \rvert + \lvert N_G(v) \cap Z \rvert + \lvert N_G(v) \cap ((Y_1 \cap V(A_{\C^*})) \cup (S^*-V(A_{\C^*} \cap B_{\C^*}))) \rvert \\
& \leq (d-1)+(s-1)+(2d-5)\\
& =  3d+s-7  = s'-1.
\end{align*}
So by the definition of $L_A$, $L_A(v) \cap \{c_B(u)\} = L_A(v) \cap L_A(u)=\emptyset$ for every $u \in Y_A \cap V(G_B) \cap N_G(v)$.
Since $N_G(v) \cap V(G_B) \subseteq Y_A$, $c_B(u) \not \in L_A(v)$ for every $u \in N_G(v) \cap V(G_B)$.
\end{proof}

Let $c$ be the $L$-coloring of $G$ defined by  $c(v):=c_A(v)$ if $v \in V(G_A)$, and $c(v):=c_B(v)$ if $v \in V(G)-V(G_A)$.

\begin{claim}
\label{no H subdivision Kst subgraph Claim13} 
Let $M$ be a monochromatic $c$-component intersecting both $V(G_A)-V(G_B)$ and $V(G_B)-V(G_A)$.
Then every component of $M \cap G_A$ intersects $Y_A$, and every component of $M \cap G_B$ intersects $Y_B$.
\end{claim}

\begin{proof}
Since $M$ intersects both $V(G_A)-V(G_B)$ and $V(G_B)-V(G_A)$, every component of $M \cap G_A$ intersects $V(A_{\C^*} \cap B_{\C^*}) \cup Z \cup S^* \cup T^* \subseteq Y_A$.

Let $M_B$ be a component of $M \cap G_B$.
Suppose that $M_B$ is disjoint from $Y_B = (Y_1 \cap V(B_{\C^*})) \cup Z \cup S^* \cup T^*$.
Since $M$ intersects both $V(G_A)-V(G_B)$ and $V(G_B)-V(G_A)$, there exist $u \in V(M_B) \cap V(A_{\C^*} \cap B_{\C^*})-S^*$ and $v \in N_{M_B}(v) \cap V(A_{\C^*}) - (V(B_{\C^*}) \cup S^* \cup T^*)$.
Since $u$ is not a $(\C^*,S^*,T^*)$-dangerous vertex, $v \not \in Y_1-T^*$.
Since $v \not \in T^*$, $v \not \in Y_1$.
So $v \in V(G_A)-(V(G_B) \cup Y_1)$ and $N_G(v) \cap V(A_{\C^*} \cap B_{\C^*})-S^* \supseteq \{u\} \neq \emptyset$.
By \cref{no H subdivision Kst subgraph Claim12}, $c(v)=c_A(v) \neq c_B(u) = c(u)$.
But $M$ is a monochromatic $c$-component, a contradiction. 
Hence every component of $M \cap G_B$ intersects $Y_B$.
\end{proof}

Let $U_A$ be the union of the monochromatic $c_A$-components of $G_A$ intersecting $Y_A$. 
Let $U_B$ be the union of the monochromatic $c_B$-components of $G_B$ intersecting $Y_B$. 
Since $c_A$ and $c_B$ are $(\eta,g)$-bounded, $\lvert V(U_A) \cup V(U_B) \rvert \leq \lvert Y_A \rvert^2g(\lvert Y_A \rvert) + \lvert Y_B \rvert^2g(\lvert Y_B \rvert) \leq (7\theta)^2g(7\theta) + (\lvert Y_1 \rvert-1)^2g(\lvert Y_1 \rvert-1) \leq g(\lvert Y_1 \rvert)$.

Since $V(G) \subseteq V(G_A) \cup V(G_B)$, by \cref{no H subdivision Kst subgraph Claim13}, every monochromatic $c$-component intersecting both $V(G_A)-V(G_B)$ and $V(G_B)-V(G_A)$ is contained in $U_A \cup U_B$ and hence contains at most $g(\lvert Y_1 \rvert) \leq \eta^2g(\eta)$ vertices.
Let $M$ be a monochromatic $c$-component. 
If $V(M) \subseteq V(G_A)$, then $M$ is a monochromatic $c_A$-component  with at most $\eta^2g(\eta)$ vertices since $c_A$ is $(\eta,g)$-bounded.
If $V(M) \subseteq V(G_B)$, then $M$ is a monochromatic $c_B$-component with at most $\eta^2g(\eta)$ vertices since $c_B$ is $(\eta,g)$-bounded.
Hence every monochromatic $c$-component contains at most $\eta^2g(\eta)$ vertices.

Since $Y_1 \subseteq Y_A \cup Y_B$, by \cref{no H subdivision Kst subgraph Claim13}, the union of the monochromatic $c$-components intersecting $Y_1$ is contained in $U_A \cup U_B$, so it contains at most $g(\lvert Y_1 \rvert) \leq \lvert Y_1 \rvert^2g(\lvert Y_1 \rvert)$ vertices. 
Therefore, $c$ is an $(\eta,g)$-bounded $L$-coloring of $G$, a contradiction.
This proves the theorem.
\end{proof}

%%%%%%%%%%%
\section{Excluding Almost $(\leq 1)$-Subdivisions}
\label{OneStepFurther}

Recall that an {\it almost $(\leq 1)$-subdivision} of a graph $H$ is a graph obtained from $H$ by subdividing edges such that at most one edge is subdivided more than once. The following simple observation is useful.

\begin{lemma} \label{cut from bipartite}
For $s\in\mathbb{N}$, let $G$ be a graph and let $H$ be a subgraph of $G$ isomorphic to $K_{s-1,t}$ for some $t \geq \binom{s-1}{2}+2$. Let $(X,Z)$ be the bipartition of $H$ with $\lvert X \rvert=s-1$. If $G$ does not contain an almost $(\leq 1)$-subdivision of $K_{s+1}$, then each component of $G-X$ contains at most one vertex in $Z$, and 
$G-X$ has at least two components.
\end{lemma}

\begin{proof}
Let $C_1,C_2,\dots,C_k$ be the components of $G-X$.
For each $i \in [k]$, $\lvert V(C_i) \cap Z \rvert \leq 1$, as otherwise $G[X \cup Z]$ together with a path in $C_i$ connecting two vertices in $V(C_i) \cap Z$ is an almost $(\leq 1)$-subdivision of $K_{s+1}$, a contradiction.
Hence $k \geq |Z| \geq t \geq 2$.
\end{proof}

The following lemma shows that a result for graphs excluding a $K_{s,t}$ subgraph can be extended for graphs 
excluding an almost $(\leq 1)$-subdivision of $K_{s+1}$.  Let $s,r \in \mathbb{N}$. Let $G$ be a graph and $Y_1 \subseteq V(G)$. 
An $(s,r,Y_1)$-list-assignment of $G$ is said to be {\it faithful} if for every $v \in V(G)-Y_1$ with $\lvert N_G(v) \cap Y_1 \rvert = s$, we have $L(v)-\bigcup_{y \in Y_1 \cap N_G(v)}L(y) \neq \emptyset$.

\begin{lemma} 
\label{push one further}
Let $\G$ be a subgraph-closed family of graphs. 
Let $\beta,r$ be functions with domain $\mathbb{N}$ such that $\beta(x) \geq x$ and $r(x) \in \mathbb{N}$ for every $x \in \mathbb{N}$.

Assume that for every $s \in \mathbb{N}$, there exist $\eta \in \mathbb{N}$ and a nondecreasing function $g$ such that for every $G \in \G$ with no $K_{s,t_s}$ subgraph, where $t_s:=\max\{ \binom{s}{2} +2, s+2\}$, 
for every $Y_1 \subseteq V(G)$ with $\lvert Y_1 \rvert \leq \eta$, and 
for every $(\beta(s),r(s),Y_1)$-list-assignment $L$ of $G$, 
there exists an $(\eta,g)$-bounded $L$-coloring of $G$.

Then for every $s \in\mathbb{N}$ with $s \geq 2$, there exist $\eta^* \in \mathbb{N}$ and a nondecreasing function $g^*$ such that 
for every graph $G \in \G$ with no almost $(\leq 1)$-subdivision of $K_{s+1}$, 
for every $Y_1 \subseteq V(G)$ with $\lvert Y_1 \rvert \leq \eta^*$, and 
for every  faithful $(\beta(s-1),r(s-1),Y_1)$-list-assignment 
$L$ of $G$, there exists an $(\eta^*,g^*)$-bounded $L$-coloring of $G$.
\end{lemma}

\begin{proof}
For every $s \in \mathbb{N}$, let $\eta_s$ be the number and $g_s$ be the function such that for every $G \in \G$ with no $K_{s,t_s}$ subgraph, every $Y_1 \subseteq V(G)$ with $\lvert Y_1 \rvert \leq \eta_s$ and every $(\beta(s),r(s),Y_1)$-list-assignment of $G$, there exists an $(\eta_s,g_s)$-bounded $L$-coloring of $G$. 
For every $s \in \mathbb{N}$ with $s \geq 2$, let $\eta_s^* := \eta_{s-1}$ and let $g_s^*$ be the function defined by $g_s^*(0) := g_{s-1}(0)$ and $g_s^*(x) := g_{s-1}(x)+\eta_{s}^* \cdot g_s^*(x-1)$ for every $x \in {\mathbb N}$. 
	
Fix $s\in\mathbb{N}-\{1\}$. Let $\beta' := \beta(s-1)$ and $r' := r(s-1)$.
We shall prove that for every graph $G$ in $\G$ with no almost $(\leq 1)$-subdivision of $K_{s+1}$, 
for every $Y_1 \subseteq V(G)$ with $\lvert Y_1 \rvert \leq \eta_s^*$, and 
for every faithful $(\beta',r',Y_1)$-list-assignment $L$ of $G$, 
there exists an $(\eta_s^*,g^*_s)$-bounded $L$-coloring of $G$.
	
Suppose to the contrary that $G$ is a graph in $\G$ with no almost $(\leq 1)$-subdivision of $K_{s+1}$, 
$Y_1$ is a subset of $V(G)$ with $\lvert Y_1 \rvert \leq \eta^*_s$, and 
$L$ is a faithful $(\beta',r',Y_1)$-list-assignment of $G$ such that 
there exists no $(\eta_s^*,g_s^*)$-bounded $L$-coloring of $G$.
We further assume that $\lvert V(G) \rvert$ is as small as possible.
	
Since $\eta_s^* = \eta_{s-1}$ and $g_s^* \geq g_{s-1}$, there exists no $(\eta_{s-1},g_{s-1})$-bounded $L$-coloring of $G$. Since $\eta_s^*=\eta_{s-1}$, by the definition of $\eta_{s-1}$ and $g_{s-1}$, $G$ contains a $K_{s-1,t_{s-1}}$ subgraph.
Let $t'$ be the maximum integer such that $G$ contains a $K_{s-1,t'}$ subgraph.
So $t' \geq t_{s-1}$.
Let $H$ be a subgraph of $G$ isomorphic to $K_{s-1,t'}$.
Let $\{P,Q\}$ be the bipartition of $H$ such that $\lvert P \rvert = s-1$ and $\lvert Q \rvert = t'$. 
By the maximality of $t'$, $Q$ is the set of all vertices of $V(G)-P$ adjacent in $G$ to all vertices in $P$.
	
	\begin{claim}
	\label{PushOneFurther1} 
	Every component of $G-P$ contains some vertex in $Y_1$.
	\end{claim}

	\begin{proof} 
	Suppose to the contrary that there exists a component $C$ of $G-P$ disjoint from $Y_1$.
	By \cref{cut from bipartite}, $G-P$ contains at least two components and there exists at most one vertex in $C$ adjacent in $G$ to all vertices in $P$.
	By the minimality of $G$, there exists an $(\eta_s^*,g_s^*)$-bounded $L|_{V(G)-V(C)}$-coloring $c$ of $G-V(C)$.
	
	Since $\beta' \geq s-1$ and $L$ is an $(\beta',r',Y_1)$-list-assignment and $V(C) \cap Y_1=\emptyset$, for every $v \in V(C)$ with $\lvert N_G(v) \cap Y_1 \rvert \leq \beta'-1$, we have
	$L(v) \cap \bigcup_{y \in N_G(v) \cap Y_1}L(y)=\emptyset$ and 
\begin{align*}
 & \;\;\;\;\; \lvert L(v)-\{c(y): y \in N_G(v) \cap P\} \rvert \\
& \geq  \lvert L(v) \rvert - \lvert \{c(y): y \in N_G(v) \cap P-Y_1\} \rvert \\
& =  \beta'+r'-\lvert N_G(v) \cap Y_1 \rvert  - \lvert \{c(y): y \in N_G(v) \cap P-Y_1\} \rvert \\
& =  \beta'+r'-\lvert N_G(v) \cap Y_1 \cap P \rvert  - \lvert \{c(y): y \in N_G(v) \cap P-Y_1\} \rvert \\
& \geq  \beta'+r'-\lvert N_G(v) \cap P \rvert \\
& \geq \beta'+r'-(s-1) \\
& \geq  1.
\end{align*}
For every $v \in V(C)$ with $\lvert N_G(v) \cap Y_1 \rvert \geq \beta'$, 
$$\lvert N_G(v) \cap P \rvert \geq \lvert N_G(v) \cap Y_1 \rvert\geq \beta' = \beta(s-1) \geq s-1,$$ 
implying $P \subseteq Y_1$ and $\lvert N_G(v) \cap Y_1 \rvert = \beta'$, so $L(v) - \{c(y): y \in N_G(v) \cap P\} = L(v) - \bigcup_{y \in N_G(v) \cap Y_1}L(y) \neq \emptyset$ since $L$ is faithful. So for every $v \in V(C)$, $L(v)-\{c(y): y \in N_G(v) \cap P\} \neq \emptyset$.
	
	Let $L'$ be the following list-assignment of $G[V(C) \cup P]$:
	\begin{itemize}
		\item For every $v \in P$, let $L'(v):=\{c(v)\}$.
		\item For every $v \in V(C)$ with $\lvert N_G(v) \cap P \rvert \geq \beta'$, let $L'(v)$ be a 1-element subset of $L(v)-\bigcup_{u \in P}L'(u)$.
		\item Let $Y_1':=P \cup \{v \in V(C): \lvert N_G(v) \cap P \rvert \geq \beta'\}$.
		\item For every $v \in V(C) \cap \nlss{Y_1'}$, let $L'(v)$ be a subset of $L(v)-\bigcup_{y \in N_G(v) \cap Y_1'}L'(y)$ of size $\lvert L(v) \rvert - \lvert N_G(v) \cap Y_1'-Y_1 \rvert = \beta'+r'-\lvert N_G(v) \cap Y_1' \rvert \geq r'+1$. 
		\item For every $v \in V(C) \cap \nlss{P}-\nlss{Y_1'}$, let $L'(v)$ be a subset of $L(v)-\bigcup_{y \in N_G(v) \cap P}L'(y)$ of size $\lvert L(v) \rvert - \lvert N_G(v) \cap P-Y_1 \rvert = \beta'+r'-\lvert N_G(v) \cap P \rvert \geq r'+1$.
		\item For every $v \in V(C)-(Y_1' \cup \nlss{Y_1'} \cup \nlss{P})$, let $L'(v):=L(v)$.
	\end{itemize}
	Note that $Y_1'-P$ consists of the vertex in $V(C)$ adjacent in $G$ to all vertices in $P$.
	Hence for every $v \in V(C) \cap N_G(P)-Y_1'$, $\lvert N_G(v) \cap P \rvert \in [ \beta'-1 ]$.
	That is, $V(C) \cap N_G(P)-Y_1' = V(C) \cap \nlss{P}$.
	So for every $v \in V(C) \cap N_G(P)-Y_1'$, $L'(v) \cap L'(u)=\emptyset$ for every $u \in P \cap N_G(v)$.
	
Clearly, $L'$ is an $(\beta',r',Y_1')$-list-assignment of $G$.
If $v \in V(C)-Y_1'$ with $\lvert N_G(v) \cap Y_1' \rvert = \beta'$, then since $\lvert Y_1'-P \rvert=1$, we know $v \in \nlss{P}-\nlss{Y_1'}$, so $L'(v)$ is a set of size at least $r'+1 \geq 2$ disjoint from $\bigcup_{y \in N_G(v) \cap P}L'(y)$.
Hence if $v \in (V(C) \cup P)-Y_1'$ with $\lvert N_G(v) \cap Y_1' \rvert = \beta'$, then $L'(v)-\bigcup_{y \in N_G(v) \cap Y_1'}L'(y) = L'(v)-\bigcup_{y \in N_G(v) \cap Y_1'-P}L'(y)$ has size $\lvert L'(y) \rvert-1 \geq 1$.
	Therefore, $L'$ is a faithful $(\beta',r',Y_1')$-list-assignment of $G[V(C) \cup P]$.
	
	Since $G-P$ contains at least two components, $\lvert V(C) \cup P \rvert < \lvert V(G) \rvert$.
	By the minimality of $G$, there exists an $(\eta_{s}^*,g_{s}^*)$-bounded $L'$-coloring $c'$ of $G[V(C) \cup P]$.
	
	For every $v \in V(C) \cap N_G(P)$, if $v \in \nlss{P}$, then $L'(v)$ is disjoint from $\bigcup_{y \in N_G(v) \cap P}L'(y)$; if $v \in V(C)$ with $\lvert N_G(v) \cap P \rvert \geq \beta'$, then $v \in Y_1'-P$ and $L'(v)$ is disjoint from $\bigcup_{y \in P}L'(y)$.
	Hence every monochromatic $c'$-component intersecting $P$ is contained in $G[P]$.
	
	Let $c^*$ be the $L$-coloring defined by $c^*(v):=c(v)$ if $v \in V(G)-V(C)$, and $c^*(v):=c'(v)$ if $v \in V(C)$.
	Hence every monochromatic $c^*$-component  is either contained in $G-V(C)$ or contained in $C$, so it contains at most ${\eta_s^*}^2g(\eta_s^*)$ vertices.
	Since $V(C) \cap Y_1=\emptyset$, the union of the monochromatic  $c^*$-components  intersecting $Y_1$ equals the union of the monochromatic  $c$-components  intersecting $Y_1$, and contains at most $\lvert Y_1 \rvert^2g(\lvert Y_1 \rvert^2)$ vertices.
	Therefore, $c^*$ is an $(\eta_s^*,g_s^*)$-bounded $L$-coloring of $G$, a contradiction.
	\end{proof}
	
	Let $C_1,C_2,\dots,C_k$ be the components of $G-P$.
	For $i \in [k]$, let $G_i:=G[V(C_i) \cup P]$.
	By \cref{cut from bipartite}, $k \geq t_{s-1} \geq s+1$. 
	By \cref{PushOneFurther1}, $Y_1 \cap V(C_i) \neq \emptyset$ for each $i \in [k]$, 
	so $\lvert V(G_i) \cap Y_1 \rvert \leq \lvert Y_1 \rvert-(k-1) \leq \lvert Y_1 \rvert-s$ for each $i \in [k]$, 
	and $k \leq \lvert Y_1 \rvert \leq \eta_s^*$.
	So $\lvert (V(G_i) \cap Y_1) \cup P \rvert \leq \lvert Y_1 \rvert -s + \lvert P \rvert < \lvert Y_1 \rvert$ for each $i \in [k]$.
	
	Let $L^*$ be the following list-assignment of $G$:
	\begin{itemize}
		\item Let $Y_1^*:=Y_1 \cup P$.
		\item For each $v \in Y_1^*$, let $L^*(v)$ be a 1-element subset of $L(v)$.
		\item For each $v \in \nlss{Y_1^*}$, let $L^*(v)$  be a subset of $L(v)-\bigcup_{y \in N_G(v) \cap Y_1^*}L^*(y)$ with size $\lvert L(v) \rvert - \lvert N_G(v) \cap (Y_1^*-Y_1) \rvert = \beta'+r'-\lvert N_G(v) \cap Y_1^* \rvert$. 
		\item For each $v \in V(G)-(Y_1^* \cup \nlss{Y_1^*})$, let $L^*(v) := L(v)$.
	\end{itemize}
	Clearly, $L^*$ is an $(\beta',r',Y_1^*)$-list-assignment of $G$.
	Let $v \in V(G)-Y_1^*$ with $\lvert N_G(v) \cap Y_1^* \rvert = \beta'$.
	So $L^*(v) = L(v)$.
If $\lvert N_G(v) \cap Y_1 \rvert = \beta'$, then $N_G(v) \cap Y_1^* = N_G(v) \cap Y_1$, so $L^*(v)-\bigcup_{y \in N_G(v) \cap Y_1^*}L^*(y) = L(v)-\bigcup_{y \in N_G(v) \cap Y_1}L(y) \neq \emptyset$ since $L$ is a faithful $(\beta',r',Y_1)$-list-assignment of $G$.
If $\lvert N_G(v) \cap Y_1 \rvert < \beta'$, then $\lvert L^*(v) \rvert = \lvert L(v) \rvert = \beta'+r'-\lvert N_G(v) \cap Y_1 \rvert$ and $L(v)$ is disjoint from $\bigcup_{y \in N_G(v) \cap Y_1}L(y)$, so 
	\begin{align*}
 \lvert L^*(v) - \bigcup_{y \in N_G(v) \cap Y_1^*}L^*(y) \rvert 
	& = \lvert L^*(v) - \bigcup_{y \in N_G(v) \cap Y_1^*-Y_1}L^*(y) \rvert\\
	& \geq \lvert L^*(v) \rvert - \lvert \bigcup_{y \in N_G(v) \cap Y_1^*-Y_1}L^*(y) \rvert \\
	& \geq \beta'+r'-\lvert N_G(v) \cap Y_1 \rvert - \lvert N_G(v) \cap Y_1^*-Y_1 \rvert\\ 
	&= \beta'+r'-\lvert N_G(v) \cap Y_1^* \rvert \\
	& = r \geq 1.
	\end{align*}
	Therefore, $L^*$ is a faithful $(\beta',r',Y_1^*)$-list-assignment of $G$.
	
	Since $P \subseteq Y_1^*$, $L^*|_{V(G_i)}$ is a faithful $(\beta',r',Y_1^* \cap V(G_i))$-list-assignment of $G_i$.
	Recall that for each $i \in [k]$, $\lvert Y_1^* \cap V(G_i) \rvert \leq \lvert Y_1 \rvert-1 \leq \eta$.
	By the minimality of $G$, for each $i \in [k]$, there exists an $(\eta_s^*,g_s^*)$-bounded $L^*|_{V(G_i)}$-coloring $c_i$ of $G_i$.
	Since $P \subseteq Y_1^* \cap V(G_i)$ for every $i \in [k]$, we know for every $v \in P$, $c_i(v)=c_j(v)$ for any $i,j \in [k]$.
	Let $c^*$ be the $L^*$-coloring of $G$ defined by $c(v):=c_1(v)$ if $v \in P$, and $c(v):=c_i(v)$ if $v \in V(C_i)$ for some $i \in [k]$.
	
	Since $P \subseteq Y_1^* \cap V(G_i)$ for all $i \in [k]$, the number of vertices in the union of the monochromatic  $c^*$-components  intersecting $Y_1 \cup P$ is at most 
	\begin{align*}
	\sum_{i=1}^k \lvert Y_1^* \cap V(G_i) \rvert^2 g_s^*(\lvert Y_1^* \cap V(G_i) \rvert) 
	& \leq \sum_{i=1}^k (\lvert Y_1 \rvert-1)^2g_s^*(\lvert Y_1 \rvert-1) \\
	& \leq \eta_s^* \cdot (\lvert Y_1 \rvert-1)^2g_s^*(\lvert Y_1 \rvert-1) \\
	& \leq \lvert Y_1 \rvert^2 g_s^*(\lvert Y_1 \rvert).
	\end{align*}
Furthermore, every monochromatic  $c^*$-component  disjoint from $Y_1 \cup P$ is a monochromatic  $c_i$-component  for some $i \in [k]$, and hence contains at most ${\eta_s^*}^2g_s^*(\eta)$ vertices.
Therefore $c^*$ is an $(\eta_s^*,g_s^*)$-bounded $L$-coloring of $G$.
This proves the lemma.
\end{proof}

The following lemma is equivalent to \cref{push one further} except it applies for restricted list assignments. The proof is identical, so we omit it. 

\begin{lemma} \label{push one further restricted}
Let $\G$ be a subgraph-closed family of graphs.
Let $\beta,r$ be functions with domain $\mathbb{N}$ such that $\beta(x) \geq x$ and $r(x) \in \mathbb{N}$ for every $x \in \mathbb{N}$.

Assume that for every $s \in \mathbb{N}$, there exist $\eta \in \mathbb{N}$ and a nondecreasing function $g$ such that for every $G \in \G$ with no $K_{s,t_s}$ subgraph, where $t_s:=\max\{ \binom{s}{2} +2, s+2\}$, 
for every $Y_1 \subseteq V(G)$ with $\lvert Y_1 \rvert \leq \eta$, and 
for every restricted $(\beta(s),r(s),Y_1)$-list-assignment $L$ of $G$, 
there exists an $(\eta,g)$-bounded $L$-coloring of $G$.

Then for every $s \in\mathbb{N}$ with $s \geq 2$, there exist $\eta^* \in \mathbb{N}$ and a nondecreasing function $g^*$ such that 
for every graph $G \in \G$ with no almost $(\leq 1)$-subdivision of $K_{s+1}$, 
for every $Y_1 \subseteq V(G)$ with $\lvert Y_1 \rvert \leq \eta^*$, and 
for every restricted faithful $(\beta(s-1),r(s-1),Y_1)$-list-assignment 
$L$ of $G$, there exists an $(\eta^*,g^*)$-bounded $L$-coloring of $G$.
\end{lemma}

\begin{theorem} \label{boosted subdiv free}
If $s\in\mathbb{N}$ with $s \geq 2$, then the following hold:
	\begin{enumerate}
		\item For every $w \in \mathbb{N}$, there exist $\eta \in \mathbb{N}$ and a nondecreasing function $g$ such that for every graph $G$ of treewidth at most $w$ with no almost $(\leq 1)$-subdivision of $K_{s+1}$, every subset $Y_1$ of $V(G)$ with $\lvert Y_1 \rvert \leq \eta$ and every faithful $(s-1,1,Y_1)$-list-assignment of $G$, there exists an $(\eta,g)$-bounded $L$-coloring of $G$.
		\item For every graph $H$, there exist $\eta \in \mathbb{N}$ and a nondecreasing function $g$ such that for every graph $G$ with no $H$-minor and no almost $(\leq 1)$-subdivision of $K_{s+1}$, every subset $Y_1$ of $V(G)$ with $\lvert Y_1 \rvert \leq \eta$ and every restricted faithful $(s-1,2,Y_1)$-list-assignment $L$ of $G$, there exists an $(\eta,g)$-bounded $L$-coloring of $G$.
		\item For every $d \in \mathbb{N}$ with $d \geq 2$ and graph $H$ of maximum degree at most $d$, there exist $\eta \in \mathbb{N}$ and a nondecreasing function $g$ such that for every graph $G$ with no $H$-subdivision and no almost $(\leq 1)$-subdivision of $K_{s+1}$, every subset $Y_1$ of $V(G)$ with $\lvert Y_1 \rvert \leq \eta$ and every restricted faithful $(s+3d-7,2,Y_1)$-list-assignment $L$ of $G$, there exists an $(\eta,g)$-bounded $L$-coloring of $G$.
		\item There exist $\eta \in \mathbb{N}$ and a nondecreasing function $g$ such that for every graph $G$ with no $K_{s+1}$-subdivision, every subset $Y_1$ of $V(G)$ with $\lvert Y_1 \rvert \leq \eta$ and every restricted faithful $(4s-7,2,Y_1)$-list-assignment $L$ of $G$, there exists an $(\eta,g)$-bounded $L$-coloring of $G$.
	\end{enumerate}
\end{theorem}

\begin{proof}
Statement 1 follows from \cref{bdd tw Kst,push one further} by taking $\G$ to be the set of graphs of treewidth at most $w$, $\beta(s)=s$ and $r(s)=1$.
	
Statement 2 follows from \cref{odd minor free,push one further restricted} by taking $\G$ to be the set of graphs with no $H$-minor, $\beta(s)=s$ and $r(s)=2$.

Statement 3 follows from \cref{push one further restricted,no H subdivision Kst subgraph} by taking $\G$ to be the set of graphs with no $H$-subdivision, $\beta(s)=3d+s-6$ and $r(s)=2$. Note that $\beta(s) \geq s$ since $d \geq 2$. And $3d+s \geq 7$ since $d \geq 2$ and $s \geq 1$.

Statement 4  follows from Statement 3 by taking $H=K_{s+1}$.
\end{proof}

When $Y_1=\emptyset$, every $(s,r,Y_1)$-list-assignment is faithful. Thus, \cref{boosted subdiv free} implies that for all $s,d,w\in\mathbb{N}$ with $s \geq 2$ and $d \geq 2$, for every graph $H$, there exists $\eta\in\mathbb{N}$ such that:
\begin{enumerate}
\item For every graph $G$ with treewidth at most $w$ and with no almost $(\leq 1)$-subdivision of $K_{s+1}$, and for every list-assignment $L$ of $G$ with $\lvert L(v) \rvert \geq s$ for every $v \in V(G)$, there exists an $L$-coloring with clustering $\eta$ (\cref{TreewidthSubdivisionChoose} for $s \geq 2$).
\item For every graph $G$ with no almost $(\leq 1)$-subdivision of $K_{s+1}$ and with no $H$-minor, there exists an $(s+1)$-coloring of $G$ with clustering $\eta$  (\cref{MinorSubdivisionColor} for $s \geq 2$).
\item If the maximum degree of $H$ is at most $d$, then for every graph $G$ with no $H$-subdivision and no almost $(\leq 1)$-subdivision of $K_{s+1}$, there exists an $(s+3d-5)$-coloring of $G$ with clustering $\eta$ (\cref{SubdivisionSubdivisionColor} for $s \geq 2$ and $d \geq 2$). 
\item For every graph $G$ with no $K_{s+1}$-subdivision, there exists a $(4s-5)$-coloring of $G$ with clustering $\eta$ (\cref{SubdivisionColor} for $s \geq 2$).
\end{enumerate}
Note that when $s=1$, graphs with no $K_{s+1}$ subgraph have no edge, so they are 1-colorable with clustering 1.
This together with \cref{easyforbidmaxdeg1subdiv} complete the proof of \cref{TreewidthSubdivisionChoose,MinorSubdivisionColor,SubdivisionSubdivisionColor,SubdivisionColor}.

%%%%%%%%%%%%%%%%%%%%%%%%%
%%%  Squashing the bibliography 
%\small
  \let\oldthebibliography=\thebibliography
  \let\endoldthebibliography=\endthebibliography
  \renewenvironment{thebibliography}[1]{%
    \begin{oldthebibliography}{#1}%
      \setlength{\parskip}{0ex}%
      \setlength{\itemsep}{0ex}%
  }{\end{oldthebibliography}}

%\bibliographystyle{DavidNatbibStyle}
%\bibliography{../../BibTeX/myBibliography.bib}

\def\soft#1{\leavevmode\setbox0=\hbox{h}\dimen7=\ht0\advance \dimen7
	by-1ex\relax\if t#1\relax\rlap{\raise.6\dimen7
		\hbox{\kern.3ex\char'47}}#1\relax\else\if T#1\relax
	\rlap{\raise.5\dimen7\hbox{\kern1.3ex\char'47}}#1\relax \else\if
	d#1\relax\rlap{\raise.5\dimen7\hbox{\kern.9ex \char'47}}#1\relax\else\if
	D#1\relax\rlap{\raise.5\dimen7 \hbox{\kern1.4ex\char'47}}#1\relax\else\if
	l#1\relax \rlap{\raise.5\dimen7\hbox{\kern.4ex\char'47}}#1\relax \else\if
	L#1\relax\rlap{\raise.5\dimen7\hbox{\kern.7ex
			\char'47}}#1\relax\else\message{accent \string\soft \space #1 not
		defined!}#1\relax\fi\fi\fi\fi\fi\fi}

\end{document}